\documentclass[final,3p,times,number]{elsarticle}
\usepackage{xcolor}
\usepackage{amsmath}
\usepackage{empheq}
\biboptions{sort&compress}
\usepackage{mathtools}
\usepackage{physics}
\usepackage{stmaryrd} % For double brackets for M jump
\usepackage{epstopdf}
\usepackage{bm}
\usepackage{graphicx,wrapfig}
\graphicspath{{./}}
\usepackage[T1]{fontenc}
\usepackage[normalem]{ulem}
\usepackage{algorithm}
\usepackage{algpseudocode}

\newcommand{\stkout}[1]{\ifmmode\text{\sout{\ensuremath{#1}}}\else\sout{#1}\fi}
\usepackage{mathrsfs}
\usepackage{amssymb}
\usepackage{nomencl}
\usepackage{etoolbox}
\usepackage{multicol}
\usepackage{caption,subcaption}
\usepackage{tabularx}
\usepackage[title]{appendix}
\usepackage{hyperref}
\newtheorem{remark}{Remark}%
\makenomenclature

\begin{document}

\begin{frontmatter}


\title{A priori error estimator for reduced-order models based on the higher-order Craig--Bampton method in structural dynamics}

\author[1]{Jaemin Kim\corref{cor1}}
\ead{jaeminkim@changwon.ac.kr}
\author[2]{Seung-Hwan Boo\corref{cor2}}
\ead{shboo@kmou.ac.kr}
\cortext[cor1]{Corresponding author}
\cortext[cor2]{Corresponding author}

\address[1]{School of Mechanical Engineering, Changwon National University, Changwon 51140, Republic of Korea}
\address[2]{Division of Naval Architecture and Ocean Systems Engineering, Korea Maritime and Ocean University, Busan 49112, Republic of Korea}

\begin{abstract}
The higher-order Craig--Bampton (HCB) method introduced by Kim et al. [AIAA J., 56 (2018)] improves reduced-model accuracy through successive residual-mode enrichments. This study uses that existing hierarchy to estimate the eigenvalue error remaining in HCB-1, a different target from the standard Craig--Bampton (CB) error addressed by earlier estimators. The proposed estimator is the relative difference between the HCB-1 and HCB-2 eigenvalues. Here, a priori denotes evaluation using both reduced eigensolutions without the full-order eigensolution. Under the stated symmetric Ritz approximation assumptions, ordered positive eigenvalues, nested augmented spaces, and preservation of retained eigenvalues by the System Equivalent Reduction Expansion Process (SEREP), the direct estimator is a non-negative lower bound on the HCB-1 relative error. This result does not require an approximation of the mode correspondence coefficient. A separate rank-1 mass-orthogonal decomposition gives a projection-based expression weighted by the squared coefficient. An exact error identity shows that close estimation requires the HCB-2 error to be small both in magnitude and relative to the HCB-1 error; the lower bound alone cannot certify an error tolerance. Extending the existing CB estimator through an HCB transformation decomposition provides a supporting result. Plate, elbow pipe, and reactor pressure vessel examples yield HCB-1 estimation-to-error ratios of approximately 0.85--0.99. The cost assessment reports the HCB-2 reference construction and solution times alongside estimator evaluation.
\end{abstract}

\end{frontmatter}

%%%%% Nomenclature Table %%%%%
\section*{Nomenclature}
\begingroup
\footnotesize
\setlength{\tabcolsep}{4pt}
\renewcommand{\arraystretch}{1.05}
\begin{center}
\begin{tabularx}{\textwidth}{l>{\raggedright\arraybackslash}X l>{\raggedright\arraybackslash}X}
\hline
Symbol & Definition & Symbol & Definition \\
\hline
$\mathbf{K}_g$, $\mathbf{M}_g$ & Global stiffness/mass matrices & $\mathbf{K}_p$, $\mathbf{M}_p$ & Partitioned (reordered) global matrices \\
$N_s$, $N_b$ & Number of substructures/boundary DOFs & $N_g$ & Total global DOFs \\
$\mathbf{K}_{ss}^{(k)}$, $\mathbf{M}_{ss}^{(k)}$ & Interior stiffness/mass of substructure $k$ & $N_d^{(k)}$, $N_t^{(k)}$ & Dominant modes/total DOFs in substructure $k$ \\
$\mathbf{K}_{sb}^{(k)}$, $\mathbf{M}_{sb}^{(k)}$ & Interior--boundary stiffness/mass coupling & $\boldsymbol{\Psi}_b^{(k)}$ & Constraint modes \\
$\boldsymbol{\Phi}_d^{(k)}$, $\boldsymbol{\Lambda}_d^{(k)}$ & Fixed-interface normal modes/eigenvalue matrix & $\hat{\mathbf{M}}_c^{(k)}$ & Coupled inertia force \\
$\mathbf{F}_i^{(k)}$ & $i$-th order residual flexibility ($i=1,2,\ldots$) & $\mathbf{R}_i^{(k)}$ & $i$-th order residual modes \\
$\mathbf{T}_{\mathrm{CB}}$, $\mathbf{T}_{\mathrm{HCB\text{-}}n}$ & CB/HCB-$n$ transformation matrices & $\mathbf{T}_{\mathrm{comp\text{-}}n}$ & Composite transformation (HCB-$n$) \\
$\mathbf{T}_0$, $\mathbf{T}_r$ & Base/residual transformation components & $\Delta\mathbf{T}_j$ & $j$-th order residual correction \\
$\mathbf{V}_{\mathrm{CB}}$, $\mathbf{V}_{\mathrm{res}}$ & SEREP matrix partition (CB/residual) & $\mathbf{T}_{\mathrm{SEREP}}$ & SEREP transformation matrix \\
$\lambda_i$, $\boldsymbol{\phi}_i$ & $i$-th eigenvalue/eigenvector (FOM) & $\tilde{\lambda}_i$, $\tilde{\boldsymbol{\phi}}_i$ & $i$-th CB eigenvalue/eigenvector \\
$\bar{\lambda}_i$, $\bar{\boldsymbol{\phi}}_i^{(n)}$ & HCB eigenpair (SEREP-reduced; reference) & $\lambda_{\mathrm{HCB},i}^{(n)}$ & $i$-th HCB-$n$ eigenvalue \\
$\xi_i$, $\xi_{\mathrm{HCB},i}^{(n)}$ & Relative eigenvalue error (CB/HCB-$n$) & $\hat{\eta}_{\mathrm{CB},i}^{(n)}$, $\hat{\eta}_{\mathrm{HCB},i}^{(1)}$ & CB/HCB-1 error estimators \\
$\mathbf{G}(\mu)$ & Auxiliary matrix $\mathbf{M}_p - \mu^{-1}\mathbf{K}_p$ & $c_i$ & Mode correspondence coefficient \\
\hline
\end{tabularx}
\end{center}
\endgroup

\section{Introduction}\label{sec:intro}
Digital twin technology has emerged as a transformative paradigm in structural engineering, enabling real-time structural health monitoring, predictive maintenance, damage detection, and remaining life prediction through continuously updated computational models \cite{torzoni2024digital, eshaghi2024realtime, kapteyn2022data}. A core requirement of digital twins is the ability to perform rapid, repeated simulations of high-fidelity finite element (FE) models, which continue to grow in size and complexity, with industrial-scale problems now routinely reaching millions of degrees of freedom (DOFs) \cite{toivanen2018multilevel, lu2021review, sun2023review}. Recent work on virtual sensing also relies on repeated FE analyses to recover full-field structural responses from sparse measurements \cite{kim2026pde}, further highlighting the need for computationally efficient structural models. However, the cost and memory requirements of full-order analyses remain prohibitive for real-time digital twin applications, driving renewed interest in reduced-order modeling (ROM) techniques. Among them, component mode synthesis (CMS)-based methods are particularly attractive because they preserve the physical substructure topology and enable independent reduction of each component, making them naturally suited for the modular architectures that underpin digital twin frameworks \cite{kapteyn2022data}. Dynamic substructuring methods have also found broad application in finite element model updating and damage detection \cite{weng2020review, zhu2021enhanced}, as well as dynamic reanalysis of structures with geometric or parametric variability \cite{mencik2023dynamic}.

Component mode synthesis (CMS) partitions a structure into substructures, reduces each independently, and assembles a reduced-order model for the global system. Among CMS techniques, the Craig--Bampton (CB) method \cite{craig1968coupling} is the most widely adopted due to its simplicity and robustness. Various extensions and accuracy studies have been reported, including iterative reduced-order substructuring approaches \cite{tian2019iterative, hagos2023modified} and comprehensive accuracy assessments of primal assembly model reduction techniques \cite{hagos2022review}. The CB method represents each substructure using fixed-interface normal modes and constraint modes. While the constraint modes capture the static response to imposed interface displacements exactly, the normal modes are truncated, introducing approximation error that can increase for higher modes. To improve accuracy, the higher-order Craig--Bampton (HCB) method \cite{kim2016development, kim2017considering} augments the CB basis with residual modes derived from the residual flexibility matrix and reduces the augmented system back to the CB size via the System Equivalent Reduction Expansion Process (SEREP) \cite{ocallahan1989serep}, achieving superior accuracy while preserving the reduced model dimension. 

An important practical question is how to assess the accuracy of a reduced model without solving the full-order problem. Without a reliable error measure, it is difficult to determine the appropriate size of the reduced model; yet if the full-order problem must be solved to evaluate the error, the computational advantage of model reduction is diminished. This fundamental paradox has long been recognized as a key obstacle to deploying reduced-order models in digital twin applications and other accuracy-critical technologies, driving sustained research efforts over the past decades. Many theoretical approaches have been developed to study eigenvalue error estimation in the context of CMS and related substructuring methods. In Ref.~\cite{elssel2006priori}, for example, an a priori bound for eigenvalue approximations obtained by the automated multilevel substructuring (AMLS) method \cite{wang2024gpu, wang2025fine, wang2026hybrid} is derived using a min--max characterization of an equivalent rational eigenvalue problem. A posteriori error estimates for CMS eigenvalue problems using residual-based error indicators have been developed \cite{jakobsson2011posteriori}, which require evaluating the eigenvalue residual with full-order matrices. For the CB method, an a priori error estimator---one that does not require solving the full-order eigenvalue problem---has been proposed \cite{kim2014estimating}, using an enhanced transformation matrix incorporating the first-order residual flexibility to predict relative eigenvalue errors. A simplified variant of this estimator was subsequently developed and applied to error control \cite{boo2016simplified}. Related developments include accuracy improvements for the dual Craig--Bampton method \cite{kim2017improving} and its extension to nonclassically damped systems \cite{gruber2018dual}.

\begin{figure}[ht!]
\centering
\includegraphics[width=0.90\textwidth]{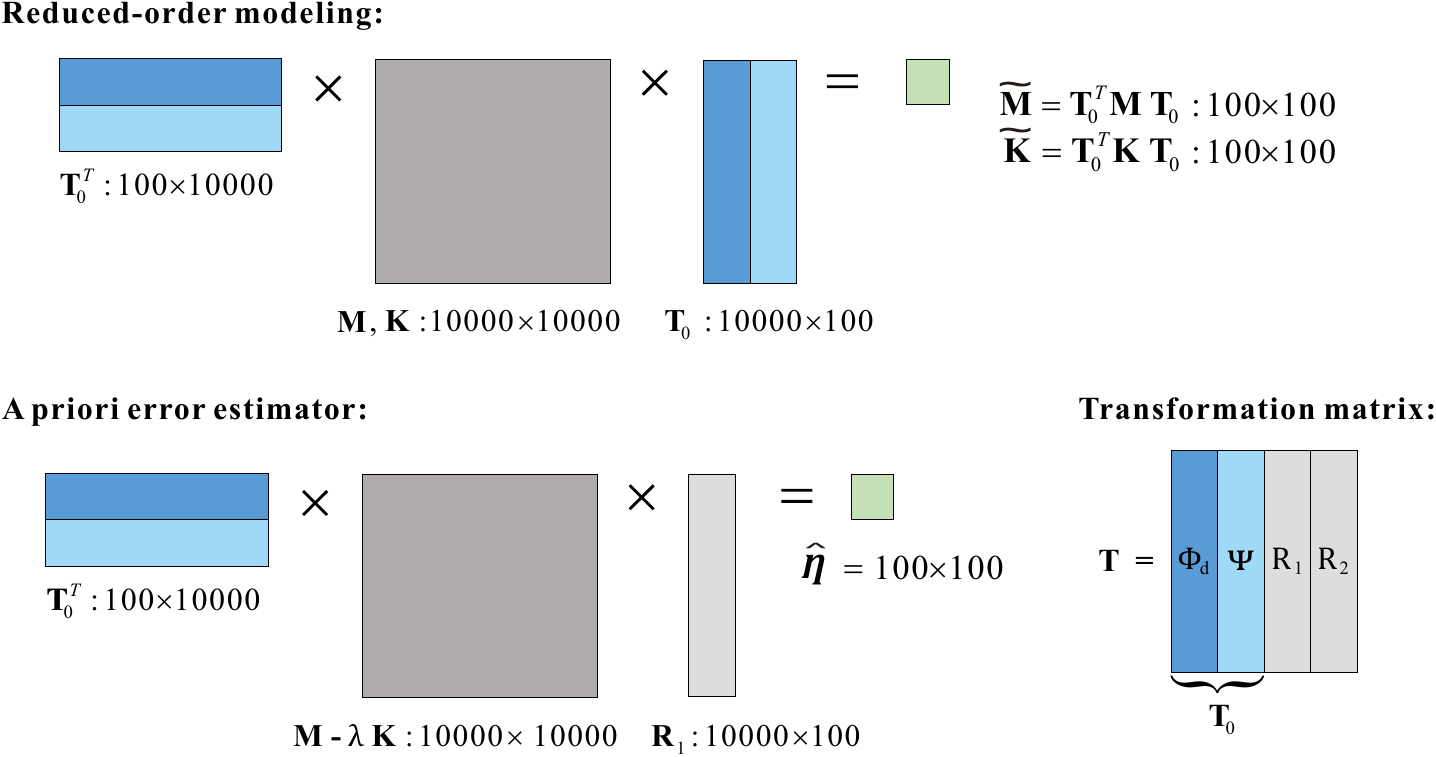}
\caption{Conceptual overview of hierarchical eigenvalue error estimation for CMS-based reduced-order models. The standard CB method is extended to the higher-order Craig--Bampton (HCB) method by augmenting the reduction basis with residual modes, and the resulting hierarchical structure is exploited to estimate eigenvalue errors of lower-order methods without solving the full-order eigenvalue problem.}
\label{fig:overview}
\end{figure}

The existing HCB method \cite{kim2017considering} supplies a hierarchy of augmented Ritz spaces by incorporating higher-order residual-mode effects. The accuracy improvement and SEREP reduction that define this hierarchy are established features of that method. The question addressed here is how to estimate the error that remains after the first-order enrichment, rather than the error of the standard CB model targeted in Refs.~\cite{kim2014estimating, boo2016simplified}. The main contribution is an HCB-1 eigenvalue-error estimator based on the HCB-2 reference, with a precise separation of its lower-bound property, its relation to a rank-1 projection expression, and the conditions under which it closely represents the actual error (Fig.~\ref{fig:overview}). The min-max bound holds for ordered positive eigenvalues of nested Ritz spaces whose retained eigenvalues are preserved by SEREP; it is independent of the approximation that the squared mode correspondence coefficient is close to unity. An exact error identity establishes the additional requirements on the HCB-2 reference error and explains why the bound alone cannot certify a tolerance. A decomposition-based extension of the CB estimator is presented as a supporting result, with its Rayleigh quotient approximations stated explicitly.

Both estimators use computed reduced eigensolutions without requiring the full-order eigensolution; the term \emph{a priori} is used in this restricted sense, consistent with the earlier CB literature. The CB formulation additionally uses full-order matrices available in the CMS setup. The computational assessment reports the HCB-2 reference construction and solution times separately from estimator evaluation. The three structural examples assess estimation accuracy at a common final ROM dimension; they do not establish a general cost or accuracy advantage over earlier estimators or an adaptive error-control procedure.

The remainder of this article is organized as follows. Section~\ref{sec:review} provides a brief review of the CB method and its higher-order extensions. Section~\ref{sec:estimator_main} develops the proposed error estimation framework. Section~\ref{sec:example} validates the proposed estimators through three numerical examples, and Section~\ref{sec:conclusion} summarizes the findings and discusses limitations. To maintain focus on the error estimation framework, supporting derivations that are not central to the main development are collected in Appendices~\ref{sec:analytical_cb}--\ref{sec:derivation_decomp}. Throughout this article, the following notational conventions are adopted. For eigenvectors, undecorated symbols ($\boldsymbol{\phi}_i$) denote exact (full-order) eigenvectors, a tilde ($\tilde{\boldsymbol{\phi}}_i$) denotes CB eigenvectors, and a bar ($\bar{\boldsymbol{\phi}}_i$) denotes SEREP-reduced HCB eigenvectors, with a superscript indicating the HCB order when necessary ($\bar{\boldsymbol{\phi}}_i^{(1)}$ for HCB-1, $\bar{\boldsymbol{\phi}}_i^{(2)}$ for HCB-2). For eigenvalues, $\lambda_i$ denotes the $i$-th exact (FOM) eigenvalue, $\tilde{\lambda}_i$ the CB eigenvalue, $\bar{\lambda}_i$ an enhanced (HCB) eigenvalue used as a reference in the error estimator, and $\lambda_{\mathrm{HCB},i}^{(n)}$ the HCB-$n$ eigenvalue. For mode-specific scalar quantities, the subscript identifies the target method and the superscript in parentheses indicates the HCB order: $\xi_i$ and $\xi_{\mathrm{HCB},i}^{(n)}$ denote the relative eigenvalue errors of the CB and HCB-$n$ methods, respectively; $\hat{\eta}_{\mathrm{CB},i}^{(n)}$ denotes the CB error estimator using HCB-$n$ as the reference; and $\hat{\eta}_{\mathrm{HCB},i}^{(1)}$ denotes the HCB-1 error estimator. A complete list of symbols is provided in the Nomenclature table.

\section{A brief review of component mode synthesis methods}\label{sec:review}

\subsection{Craig--Bampton (CB) method}\label{sec:cb}

The undamped free vibration of a structure is governed by the equations of motion
\begin{align}\label{eq:eom}
    \mathbf{M}_g \ddot{\mathbf{u}} + \mathbf{K}_g \mathbf{u} = \mathbf{0},
\end{align}
where $\mathbf{M}_g$ and $\mathbf{K}_g$ are the global mass and stiffness matrices, and $\mathbf{u}$ is the displacement vector. Assuming harmonic motion $\mathbf{u} = \boldsymbol{\phi} e^{i\omega t}$ with natural frequency $\omega$, Eq.~\eqref{eq:eom} reduces to
\begin{align}\label{eq:eigproblem}
    \mathbf{K}_g \boldsymbol{\phi} = \lambda \mathbf{M}_g \boldsymbol{\phi},
\end{align}
where $\lambda = \omega^2$ is the eigenvalue and $\boldsymbol{\phi}$ is the corresponding eigenvector. For large-scale models, solving Eq.~\eqref{eq:eigproblem} directly is prohibitively expensive, motivating the use of CMS reduction.

The structure is divided into $N_s$ substructures. The last substructure ($k = N_s$) collects all boundary (interface) degrees of freedom (DOFs), while the remaining substructures ($k = 1, \ldots, N_s - 1$) represent interior regions. For each interior substructure $k$, the global matrices are partitioned as
\begin{align}\label{eq:partition}
    \mathbf{K}^{(k)} = \begin{bmatrix} \mathbf{K}_{ss}^{(k)} & \mathbf{K}_{sb}^{(k)} \\ \mathbf{K}_{bs}^{(k)} & \mathbf{K}_{bb}^{(k)} \end{bmatrix}, \quad
    \mathbf{M}^{(k)} = \begin{bmatrix} \mathbf{M}_{ss}^{(k)} & \mathbf{M}_{sb}^{(k)} \\ \mathbf{M}_{bs}^{(k)} & \mathbf{M}_{bb}^{(k)} \end{bmatrix}, \quad
    \mathbf{u}^{(k)} = \begin{bmatrix} \mathbf{u}_s^{(k)} \\ \mathbf{u}_b \end{bmatrix}, \quad
    \ddot{\mathbf{u}}^{(k)} = \begin{bmatrix} \ddot{\mathbf{u}}_s^{(k)} \\ \ddot{\mathbf{u}}_b \end{bmatrix},
\end{align}
where subscript $s$ denotes interior DOFs of substructure $k$ and subscript $b$ denotes boundary DOFs. After reordering, the global system takes the block form
\begin{align}\label{eq:blockassembly}
    \mathbf{K}_p = \begin{bmatrix}
    \mathbf{K}_{ss}^{(1)} & & & \mathbf{K}_{sb}^{(1)} \\
    & \ddots & & \vdots \\
    & & \mathbf{K}_{ss}^{(N_s-1)} & \mathbf{K}_{sb}^{(N_s-1)} \\
    \mathbf{K}_{bs}^{(1)} & \cdots & \mathbf{K}_{bs}^{(N_s-1)} & \mathbf{K}_{bb}
    \end{bmatrix}, \quad
    \mathbf{u}_p = \begin{bmatrix} \mathbf{u}_s^{(1)} \\ \vdots \\ \mathbf{u}_s^{(N_s-1)} \\ \mathbf{u}_b \end{bmatrix}, \quad
    \ddot{\mathbf{u}}_p = \begin{bmatrix} \ddot{\mathbf{u}}_s^{(1)} \\ \vdots \\ \ddot{\mathbf{u}}_s^{(N_s-1)} \\ \ddot{\mathbf{u}}_b \end{bmatrix},
\end{align}
and similarly for $\mathbf{M}_p$, where $\mathbf{K}_{bb} = \sum_{k} \mathbf{K}_{bb}^{(k)}$ and $\mathbf{M}_{bb} = \sum_{k} \mathbf{M}_{bb}^{(k)}$.
From Eq.~\eqref{eq:eom}, the equation of motion for the interior DOFs of substructure $k$ is
\begin{align}\label{eq:eom_sub}
    \mathbf{M}_{ss}^{(k)} \ddot{\mathbf{u}}_s^{(k)} + \mathbf{M}_{sb}^{(k)} \ddot{\mathbf{u}}_b + \mathbf{K}_{ss}^{(k)} \mathbf{u}_s^{(k)} + \mathbf{K}_{sb}^{(k)} \mathbf{u}_b = \mathbf{0}.
\end{align}

For each interior substructure $k$, the fixed-interface eigenvalue problem is
\begin{align}\label{eq:fixed_eig}
    \mathbf{K}_{ss}^{(k)} \boldsymbol{\phi}_j^{(k)} = \lambda_j^{(k)} \mathbf{M}_{ss}^{(k)} \boldsymbol{\phi}_j^{(k)}, \qquad j = 1, \ldots, N_d^{(k)},
\end{align}
where $N_d^{(k)}$ is the number of dominant modes, selected in ascending order of eigenvalue. The lower fixed-interface modes are retained to represent the low-frequency response with a limited number of basis vectors. The eigenvectors are collected into $\boldsymbol{\Phi}_d^{(k)} = [\boldsymbol{\phi}_1^{(k)}, \ldots, \boldsymbol{\phi}_{N_d^{(k)}}^{(k)}]$ and mass-normalized such that
\begin{align}\label{eq:massnorm}
    (\boldsymbol{\Phi}_d^{(k)})^{\mathrm{T}} \mathbf{M}_{ss}^{(k)} \boldsymbol{\Phi}_d^{(k)} = \mathbf{I}, \qquad
    (\boldsymbol{\Phi}_d^{(k)})^{\mathrm{T}} \mathbf{K}_{ss}^{(k)} \boldsymbol{\Phi}_d^{(k)} = \boldsymbol{\Lambda}_d^{(k)} = \mathrm{diag}(\lambda_1^{(k)}, \ldots, \lambda_{N_d^{(k)}}^{(k)}).
\end{align}

The constraint modes represent the static response of interior DOFs to unit displacements at boundary DOFs while all other boundaries are fixed:
\begin{align}\label{eq:constraint}
    \boldsymbol{\Psi}_b^{(k)} = -\left(\mathbf{K}_{ss}^{(k)}\right)^{-1} \mathbf{K}_{sb}^{(k)}.
\end{align}

The CB transformation for substructure $k$ relates the physical DOFs to the generalized coordinates $(\mathbf{q}_d^{(k)}, \mathbf{u}_b)$:
\begin{align}\label{eq:cb_trans_sub}
    \begin{bmatrix} \mathbf{u}_s^{(k)} \\ \mathbf{u}_b \end{bmatrix}
    = \begin{bmatrix} \boldsymbol{\Phi}_d^{(k)} & \boldsymbol{\Psi}_b^{(k)} \\ \mathbf{0} & \mathbf{I} \end{bmatrix}
    \begin{bmatrix} \mathbf{q}_d^{(k)} \\ \mathbf{u}_b \end{bmatrix}.
\end{align}
Assembling over all substructures yields the global CB transformation matrix $\mathbf{T}_{\mathrm{CB}}$:
\begin{align}\label{eq:cb_trans_global}
    \mathbf{T}_{\mathrm{CB}} = \begin{bmatrix}
    \boldsymbol{\Phi}_d^{(1)} & & & \boldsymbol{\Psi}_b^{(1)} \\
    & \ddots & & \vdots \\
    & & \boldsymbol{\Phi}_d^{(N_s-1)} & \boldsymbol{\Psi}_b^{(N_s-1)} \\
    \mathbf{0} & \cdots & \mathbf{0} & \mathbf{I}
    \end{bmatrix}.
\end{align}
The CB reduced matrices are obtained by the congruence transformation
\begin{align}\label{eq:cb_reduced}
    \widetilde{\mathbf{K}}_{\mathrm{CB}} = \mathbf{T}_{\mathrm{CB}}^{\mathrm{T}} \mathbf{K}_p \mathbf{T}_{\mathrm{CB}}, \qquad
    \widetilde{\mathbf{M}}_{\mathrm{CB}} = \mathbf{T}_{\mathrm{CB}}^{\mathrm{T}} \mathbf{M}_p \mathbf{T}_{\mathrm{CB}},
\end{align}
yielding a reduced eigenvalue problem of size $N_d + N_b$, where $N_d = \sum_{k=1}^{N_s-1} N_d^{(k)}$ is the total number of dominant modes and $N_b$ is the number of boundary DOFs.

\subsection{Higher-order Craig--Bampton (HCB) method}\label{sec:hcb}
The key idea of the HCB method is to explicitly augment the CB transformation basis with residual modes derived from the residual flexibility. Since the augmented basis introduces additional generalized coordinates (one set of $N_b$ coordinates per residual mode order), the augmented system is larger than the original CB system. To maintain the same reduced-model dimension $N_d + N_b$ as the CB method---and thereby enable a fair comparison at equal model size---the augmented system is projected back to the CB size via SEREP (System Equivalent Reduction Expansion Process).

\subsubsection{Higher-order residual flexibility}\label{sec:hcb_derivation}
To motivate the residual modes, consider the full set of fixed-interface normal modes partitioned into $N_d^{(k)}$ dominant modes $\boldsymbol{\Phi}_d^{(k)}$ and $N_r^{(k)} = N_t^{(k)} - N_d^{(k)}$ residual modes $\boldsymbol{\Phi}_r^{(k)}$. Including both, the full transformation relates the physical interior displacement to the modal coordinates:
\begin{align}\label{eq:full_interior}
    \mathbf{u}_s^{(k)} = \boldsymbol{\Phi}_d^{(k)} \mathbf{q}_d^{(k)} + \boldsymbol{\Phi}_r^{(k)} \mathbf{q}_r^{(k)} + \boldsymbol{\Psi}_b^{(k)} \mathbf{u}_b.
\end{align}
Substituting Eq.~\eqref{eq:full_interior} into the equation of motion~\eqref{eq:eom_sub}, pre-multiplying by $(\boldsymbol{\Phi}_r^{(k)})^{\mathrm{T}}$, and assuming harmonic response ($\lambda = \omega^2$), the residual modal coordinates are obtained as (see Appendix~\ref{sec:derivation_qr} for details)
\begin{align}\label{eq:qr_solution}
    \mathbf{q}_r^{(k)} = \lambda \left[\boldsymbol{\Lambda}_r^{(k)} - \lambda \mathbf{I}\right]^{-1} (\boldsymbol{\Phi}_r^{(k)})^{\mathrm{T}} \hat{\mathbf{M}}_c^{(k)} \mathbf{u}_b,
\end{align}
where the coupled inertia force is defined as
\begin{align}\label{eq:Mc_hat}
    \hat{\mathbf{M}}_c^{(k)} \coloneqq \mathbf{M}_{ss}^{(k)} \boldsymbol{\Psi}_b^{(k)} + \mathbf{M}_{sb}^{(k)},
\end{align}
representing the inertia force acting on the interior DOFs due to boundary acceleration. Substituting Eq.~\eqref{eq:qr_solution} into Eq.~\eqref{eq:full_interior}, the correction to the constraint mode involves the residual flexibility:
\begin{align}\label{eq:us_corrected}
    \mathbf{u}_s^{(k)} = \boldsymbol{\Phi}_d^{(k)} \mathbf{q}_d^{(k)} + \boldsymbol{\Psi}_b^{(k)} \mathbf{u}_b + \lambda \boldsymbol{\Phi}_r^{(k)} \left(\boldsymbol{\Lambda}_r^{(k)} - \lambda \mathbf{I}\right)^{-1} (\boldsymbol{\Phi}_r^{(k)})^{\mathrm{T}} \hat{\mathbf{M}}_c^{(k)} \mathbf{u}_b.
\end{align}
The residual flexibility term can be expanded via a Neumann series:
\begin{align}\label{eq:taylor_expansion}
    \boldsymbol{\Phi}_r^{(k)} \left[\boldsymbol{\Lambda}_r^{(k)} - \lambda \mathbf{I}\right]^{-1} (\boldsymbol{\Phi}_r^{(k)})^{\mathrm{T}}
    &= \underbrace{\boldsymbol{\Phi}_r^{(k)} \left(\boldsymbol{\Lambda}_r^{(k)}\right)^{-1} (\boldsymbol{\Phi}_r^{(k)})^{\mathrm{T}}}_{\coloneqq\, \mathbf{F}_1^{(k)}}
     + \lambda \underbrace{\boldsymbol{\Phi}_r^{(k)} \left(\boldsymbol{\Lambda}_r^{(k)}\right)^{-2} (\boldsymbol{\Phi}_r^{(k)})^{\mathrm{T}}}_{\coloneqq\, \mathbf{F}_2^{(k)}}
     + \cdots
\end{align}
where
\begin{align}\label{eq:F_i_def}
    \mathbf{F}_i^{(k)} = \boldsymbol{\Phi}_r^{(k)} \left(\boldsymbol{\Lambda}_r^{(k)}\right)^{-i} (\boldsymbol{\Phi}_r^{(k)})^{\mathrm{T}}, \qquad i = 1, 2, \ldots
\end{align}
is the $i$-th order residual flexibility (see Remark~\ref{rmk:neumann} for the Neumann series derivation and convergence). In practice, $\mathbf{F}_1^{(k)}$ can be computed without the residual eigenvalues and eigenvectors as
\begin{align}\label{eq:res_flex}
    \mathbf{F}_1^{(k)} = \left(\mathbf{K}_{ss}^{(k)}\right)^{-1} - \boldsymbol{\Phi}_d^{(k)} \left(\boldsymbol{\Lambda}_d^{(k)}\right)^{-1} (\boldsymbol{\Phi}_d^{(k)})^{\mathrm{T}}.
\end{align}
This identity follows from the spectral decomposition of the stiffness inverse (see Appendix~\ref{sec:derivation_Frs} for details). Since the complete eigenvectors $[\boldsymbol{\Phi}_d^{(k)}, \boldsymbol{\Phi}_r^{(k)}]$ diagonalize $\mathbf{K}_{ss}^{(k)}$, the full flexibility is
\begin{align*}
    (\mathbf{K}_{ss}^{(k)})^{-1} = \boldsymbol{\Phi}_d^{(k)} (\boldsymbol{\Lambda}_d^{(k)})^{-1} (\boldsymbol{\Phi}_d^{(k)})^{\mathrm{T}} + \boldsymbol{\Phi}_r^{(k)} (\boldsymbol{\Lambda}_r^{(k)})^{-1} (\boldsymbol{\Phi}_r^{(k)})^{\mathrm{T}},
\end{align*}
i.e., the sum of the dominant and residual contributions. Subtracting the known dominant part isolates the residual flexibility $\mathbf{F}_1^{(k)}$ without requiring the residual eigenpairs explicitly.

The second-order residual flexibility $\mathbf{F}_2^{(k)} = \boldsymbol{\Phi}_r^{(k)} (\boldsymbol{\Lambda}_r^{(k)})^{-2} (\boldsymbol{\Phi}_r^{(k)})^{\mathrm{T}}$ (Eq.~\eqref{eq:F_i_def} with $i = 2$) can also be computed from the first-order (Eq.~\eqref{eq:F_i_def} with $i = 1$) as
\begin{align}\label{eq:Frm}
    \mathbf{F}_2^{(k)} = \mathbf{F}_1^{(k)} \mathbf{M}_{ss}^{(k)} \mathbf{F}_1^{(k)}
    = \boldsymbol{\Phi}_r^{(k)} \bigl(\boldsymbol{\Lambda}_r^{(k)}\bigr)^{-1}
       \underbrace{(\boldsymbol{\Phi}_r^{(k)})^{\mathrm{T}} \mathbf{M}_{ss}^{(k)} \boldsymbol{\Phi}_r^{(k)}}_{= \, \mathbf{I}}
       \bigl(\boldsymbol{\Lambda}_r^{(k)}\bigr)^{-1} (\boldsymbol{\Phi}_r^{(k)})^{\mathrm{T}}
    = \boldsymbol{\Phi}_r^{(k)} \bigl(\boldsymbol{\Lambda}_r^{(k)}\bigr)^{-2} (\boldsymbol{\Phi}_r^{(k)})^{\mathrm{T}},
\end{align}
where the mass-orthonormality of the residual eigenvectors, $(\boldsymbol{\Phi}_r^{(k)})^{\mathrm{T}} \mathbf{M}_{ss}^{(k)} \boldsymbol{\Phi}_r^{(k)} = \mathbf{I}$, eliminates the central factor.

\begin{remark}[Neumann series expansion and convergence]\label{rmk:neumann}\label{rmk:convergence}
    The Neumann series states that
    \begin{align}
        (\mathbf{I} - \mathbf{B})^{-1} = \sum_{j=0}^{\infty} \mathbf{B}^j, \qquad \|\mathbf{B}\| < 1. \notag
    \end{align}
    In Eq.~\eqref{eq:taylor_expansion}, this is applied with
    \begin{align}
        [\boldsymbol{\Lambda}_r^{(k)} - \lambda \mathbf{I}]^{-1} = \left(\boldsymbol{\Lambda}_r^{(k)}\right)^{-1} \left[\mathbf{I} - \lambda \left(\boldsymbol{\Lambda}_r^{(k)}\right)^{-1}\right]^{-1}, \qquad \mathbf{B} = \lambda \left(\boldsymbol{\Lambda}_r^{(k)}\right)^{-1}, \notag
    \end{align}
    yielding an expansion of the dynamic residual flexibility in ascending powers of $\lambda$, where the first-order term is $\mathbf{F}_1^{(k)}$ and the second-order term is $\mathbf{F}_2^{(k)}$, each capturing progressively higher-order inertia effects of the truncated modes. The present work limits the expansion to order $n$, defining the HCB-$n$ method; hence the zeroth-order truncation (no residual correction) recovers the standard CB method, the first-order truncation yields HCB-1, and so forth. This systematic hierarchy is central to the proposed error estimation framework: the $(n+1)$-th order solution provides a more complete approximation that can serve as a reference for assessing the $n$-th order error.

    The series converges if $\lambda < \min(\mathrm{diag}(\boldsymbol{\Lambda}_r^{(k)}))$, i.e., the eigenvalue of interest is smaller than the smallest residual eigenvalue. The applicability of the expansion therefore depends on this condition holding for each substructure over the eigenvalue range of interest.
\end{remark}

\subsubsection{Residual modes}\label{sec:residual_modes}
Truncating the Neumann series in Eq.~\eqref{eq:us_corrected} at order $n$, the corrected interior displacement becomes
\begin{align}\label{eq:us_truncated}
    \mathbf{u}_s^{(k)} \approx \boldsymbol{\Phi}_d^{(k)} \mathbf{q}_d^{(k)} + \boldsymbol{\Psi}_b^{(k)} \mathbf{u}_b + \sum_{i=1}^{n} \lambda^i \mathbf{F}_i^{(k)} \hat{\mathbf{M}}_c^{(k)} \mathbf{u}_b.
\end{align}
By introducing the generalized coordinates $\boldsymbol{\eta}_i = \lambda^i \mathbf{u}_b$ ($i = 1, \ldots, n$), each correction term $\mathbf{F}_i^{(k)} \hat{\mathbf{M}}_c^{(k)} \boldsymbol{\eta}_i$ takes the form of a column-space projection analogous to the constraint mode $\boldsymbol{\Psi}_b^{(k)} \mathbf{u}_b$. The $i$-th order residual mode is thus defined as
\begin{align}\label{eq:residual_mode_general}
    \mathbf{R}_i^{(k)} \coloneqq \mathbf{F}_i^{(k)} \hat{\mathbf{M}}_c^{(k)}, \qquad i = 1, \ldots, n,
\end{align}
which represents the deformation pattern of the truncated modes under $i$-th order inertia loading filtered through $\mathbf{F}_i^{(k)}$. The first-order residual mode $\mathbf{R}_1^{(k)}$ captures the quasi-static deformation of the truncated modes under the inertia loading, and the second-order residual mode $\mathbf{R}_2^{(k)} = \mathbf{F}_2^{(k)} \hat{\mathbf{M}}_c^{(k)}$ accounts for the next-order inertia effect by filtering the force through the mass-weighted residual flexibility.

The columns of each residual mode matrix are $L_2$-normalized. Denoting the $j$-th column of $\mathbf{R}_i^{(k)}$ by $\mathbf{r}_j$, the normalized column is
\begin{align}
    \hat{\mathbf{r}}_j = \frac{\mathbf{r}_j}{\|\mathbf{r}_j\|_2}, \qquad j = 1, \ldots, N_b.
\end{align}
The columns of $\mathbf{F}_i^{(k)} \hat{\mathbf{M}}_c^{(k)}$ can differ by several orders of magnitude depending on the boundary DOF orientation and the local stiffness. Large disparities in column norms can increase the condition number $\kappa(\widetilde{\mathbf{K}}) = \|\widetilde{\mathbf{K}}\| \, \|\widetilde{\mathbf{K}}^{-1}\|$ of the augmented system matrices, amplifying rounding errors in the eigenvalue solve. Normalizing the columns to unit length reduces these differences in scale.

The $n$-th order HCB method (HCB-$n$) augments the CB basis with residual modes $\mathbf{R}_1^{(k)}, \ldots, \mathbf{R}_n^{(k)}$. The corrected interior displacement in Eq.~\eqref{eq:us_truncated} is therefore expressed as
\begin{align}\label{eq:us_hcb_n}
    \mathbf{u}_s^{(k)} \approx \boldsymbol{\Phi}_d^{(k)} \mathbf{q}_d^{(k)} + \boldsymbol{\Psi}_b^{(k)} \mathbf{u}_b + \sum_{i=1}^{n} \mathbf{R}_i^{(k)} \boldsymbol{\eta}_i.
\end{align}
HCB-0 (i.e., $n = 0$) reduces to the standard CB method.

\subsubsection{HCB-1: First-order residual mode augmentation}\label{sec:hcb1}
The HCB-1 transformation augments the CB basis with $\mathbf{R}_1^{(k)}$:
\begin{align}\label{eq:hcb1_trans_sub}
    \begin{bmatrix} \mathbf{u}_s^{(k)} \\ \mathbf{u}_b \end{bmatrix}
    = \begin{bmatrix} \boldsymbol{\Phi}_d^{(k)} & \boldsymbol{\Psi}_b^{(k)} & \mathbf{R}_1^{(k)} \\ \mathbf{0} & \mathbf{I} & \mathbf{0} \end{bmatrix}
    \begin{bmatrix} \mathbf{q}_d^{(k)} \\ \mathbf{u}_b \\ \boldsymbol{\eta}_1 \end{bmatrix},
\end{align}
where $\boldsymbol{\eta}_1 = \lambda \mathbf{u}_b$ are the generalized coordinates associated with the residual mode. The global transformation is
\begin{align}\label{eq:hcb1_trans_global}
    \mathbf{T}_{\mathrm{HCB\text{-}1}} = \begin{bmatrix}
    \boldsymbol{\Phi}_d^{(1)} & & & \boldsymbol{\Psi}_b^{(1)} & \mathbf{R}_1^{(1)} \\
    & \ddots & & \vdots & \vdots \\
    & & \boldsymbol{\Phi}_d^{(N_s-1)} & \boldsymbol{\Psi}_b^{(N_s-1)} & \mathbf{R}_1^{(N_s-1)} \\
    \mathbf{0} & \cdots & \mathbf{0} & \mathbf{I} & \mathbf{0}
    \end{bmatrix}
    = \bigl[\, \mathbf{T}_{\mathrm{CB}} \;\big|\; \Delta\mathbf{T}_1 \,\bigr],
\end{align}
where $\Delta\mathbf{T}_1$ denotes the assembled first-order residual correction whose $k$-th block row is $[\mathbf{0}, \ldots, \mathbf{0}, \mathbf{R}_1^{(k)}]$. This notation explicitly separates the CB basis from the residual augmentation and will be reused in Section~\ref{sec:estimator} for constructing the a priori error estimator.

The reduced matrices are assembled in the block form
\begin{align}\label{eq:hcb1_reduced}
    \widetilde{\mathbf{K}}_{\mathrm{HCB\text{-}1}} = \begin{bmatrix} \widetilde{\mathbf{K}}_{\mathrm{CB}} & \mathbf{K}_{c,r_1} \\ \mathbf{K}_{r_1,c} & \mathbf{K}_{r_1,r_1} \end{bmatrix}, \qquad
    \widetilde{\mathbf{M}}_{\mathrm{HCB\text{-}1}} = \begin{bmatrix} \widetilde{\mathbf{M}}_{\mathrm{CB}} & \mathbf{M}_{c,r_1} \\ \mathbf{M}_{r_1,c} & \mathbf{M}_{r_1,r_1} \end{bmatrix},
\end{align}
where the coupling blocks are
\begin{align}
    \mathbf{K}_{r_1,r_1} = \sum_{k=1}^{N_s-1} (\mathbf{R}_1^{(k)})^{\mathrm{T}} \mathbf{K}_{ss}^{(k)} \mathbf{R}_1^{(k)}, \qquad
    \mathbf{M}_{r_1,r_1} = \sum_{k=1}^{N_s-1} (\mathbf{R}_1^{(k)})^{\mathrm{T}} \mathbf{M}_{ss}^{(k)} \mathbf{R}_1^{(k)},
\end{align}
and the off-diagonal blocks involve coupling terms $(\boldsymbol{\Phi}_d^{(k)})^{\mathrm{T}} \mathbf{K}_{ss}^{(k)} \mathbf{R}_1^{(k)}$, $((\boldsymbol{\Psi}_b^{(k)})^{\mathrm{T}} \mathbf{K}_{ss}^{(k)} + \mathbf{K}_{bs}^{(k)}) \mathbf{R}_1^{(k)}$, and analogous mass terms.
The augmented system has size $\sum_k N_d^{(k)} + 2N_b$, which is larger than the CB system.

\subsubsection{HCB-2: First- and second-order residual mode augmentation}\label{sec:hcb2}
The HCB-2 transformation further includes $\mathbf{R}_2^{(k)}$:
\begin{align}\label{eq:hcb2_trans_sub}
    \begin{bmatrix} \mathbf{u}_s^{(k)} \\ \mathbf{u}_b \end{bmatrix}
    = \begin{bmatrix} \boldsymbol{\Phi}_d^{(k)} & \boldsymbol{\Psi}_b^{(k)} & \mathbf{R}_1^{(k)} & \mathbf{R}_2^{(k)} \\ \mathbf{0} & \mathbf{I} & \mathbf{0} & \mathbf{0} \end{bmatrix}
    \begin{bmatrix} \mathbf{q}_d^{(k)} \\ \mathbf{u}_b \\ \boldsymbol{\eta}_1 \\ \boldsymbol{\eta}_2 \end{bmatrix},
\end{align}
where $\boldsymbol{\eta}_2 = \lambda^2 \mathbf{u}_b$ is the generalized coordinate associated with the second-order residual mode. The global transformation is
\begin{align}\label{eq:hcb2_trans_global}
    \mathbf{T}_{\mathrm{HCB\text{-}2}} = \begin{bmatrix}
    \boldsymbol{\Phi}_d^{(1)} & & & \boldsymbol{\Psi}_b^{(1)} & \mathbf{R}_1^{(1)} & \mathbf{R}_2^{(1)} \\
    & \ddots & & \vdots & \vdots & \vdots \\
    & & \boldsymbol{\Phi}_d^{(N_s-1)} & \boldsymbol{\Psi}_b^{(N_s-1)} & \mathbf{R}_1^{(N_s-1)} & \mathbf{R}_2^{(N_s-1)} \\
    \mathbf{0} & \cdots & \mathbf{0} & \mathbf{I} & \mathbf{0} & \mathbf{0}
    \end{bmatrix}
    = \bigl[\, \mathbf{T}_{\mathrm{CB}} \;\big|\; \Delta\mathbf{T}_1 \;\big|\; \Delta\mathbf{T}_2 \,\bigr],
\end{align}
where $\Delta\mathbf{T}_2$ denotes the assembled second-order residual correction whose $k$-th block row is $[\mathbf{0}, \ldots, \mathbf{0}, \mathbf{R}_2^{(k)}]$. The reduced matrices are assembled in the extended block form
\begin{align}\label{eq:hcb2_reduced}
    \widetilde{\mathbf{K}}_{\mathrm{HCB\text{-}2}} = \begin{bmatrix}
    \widetilde{\mathbf{K}}_{\mathrm{CB}} & \mathbf{K}_{c,r_1} & \mathbf{K}_{c,r_2} \\
    \mathbf{K}_{r_1,c} & \mathbf{K}_{r_1,r_1} & \mathbf{K}_{r_1,r_2} \\
    \mathbf{K}_{r_2,c} & \mathbf{K}_{r_2,r_1} & \mathbf{K}_{r_2,r_2}
    \end{bmatrix},
\end{align}
and similarly for $\widetilde{\mathbf{M}}_{\mathrm{HCB\text{-}2}}$. The augmented system has size $\sum_k N_d^{(k)} + 3N_b$.

\subsubsection{SEREP reduction}\label{sec:serep}
The augmented transformation $\mathbf{T}_{\mathrm{HCB\text{-}}n} = [\, \mathbf{T}_{\mathrm{CB}} \mid \Delta\mathbf{T}_1 \mid \cdots \mid \Delta\mathbf{T}_n \,]$ has $n_{\mathrm{CB}} + n N_b$ columns, which is larger than the CB system size $n_{\mathrm{CB}}$ due to the residual mode coordinates $\boldsymbol{\eta}_1, \ldots, \boldsymbol{\eta}_n$. To restore the reduced-model size to $n_{\mathrm{CB}}$, a secondary reduction based on the concept of the System Equivalent Reduction Expansion Process (SEREP) \cite{ocallahan1989serep} is applied \cite{kim2017considering}. The original SEREP reduces a full-order model to a subset of physical DOFs using measured or analytical mode shapes. Here, it reduces the augmented HCB-$n$ system to the CB dimension within the reduced-order space. This secondary reduction also enables a fair comparison with the CB method at the same model dimension. The SEREP transformation is constructed by solving the eigenvalue problem of the augmented system and retaining the first $n_{\mathrm{CB}}$ eigenpairs:
\begin{align}\label{eq:serep_eig}
    \widetilde{\mathbf{K}}_{\mathrm{HCB\text{-}}n} \boldsymbol{\psi}_j = \mu_j \widetilde{\mathbf{M}}_{\mathrm{HCB\text{-}}n} \boldsymbol{\psi}_j, \qquad j = 1, \ldots, n_{\mathrm{CB}}, \quad n = 1, 2, \ldots
\end{align}
where $n_{\mathrm{CB}} = \sum_k N_d^{(k)} + N_b$ is the CB reduced-model size. Collecting the retained eigenvectors into $\mathbf{T}_{\mathrm{SEREP}} = [\boldsymbol{\psi}_1, \ldots, \boldsymbol{\psi}_{n_{\mathrm{CB}}}]$, the final reduced matrices are
\begin{align}\label{eq:serep_reduced}
    \overline{\mathbf{K}} = \mathbf{T}_{\mathrm{SEREP}}^{\mathrm{T}} \widetilde{\mathbf{K}}_{\mathrm{HCB\text{-}}n} \mathbf{T}_{\mathrm{SEREP}}, \qquad
    \overline{\mathbf{M}} = \mathbf{T}_{\mathrm{SEREP}}^{\mathrm{T}} \widetilde{\mathbf{M}}_{\mathrm{HCB\text{-}}n} \mathbf{T}_{\mathrm{SEREP}}.
\end{align}
The final eigenvalue problem $\overline{\mathbf{K}} \bar{\boldsymbol{\phi}} = \bar{\lambda} \overline{\mathbf{M}} \bar{\boldsymbol{\phi}}$ has the same size as the CB problem, but with improved accuracy due to the residual mode information embedded through SEREP.

\section{Hierarchical eigenvalue error estimation for HCB models}\label{sec:estimator_main}
This section develops the HCB-1 estimator using HCB-2 as a reference, preceded by a supporting extension of the CB estimator in \cite{kim2014estimating}. Both formulations use computed reduced eigensolutions without requiring the full-order eigensolution. Here, the term \emph{a priori} used in the earlier CB literature and in Fig.~\ref{fig:overview} refers only to this independence from the full-order eigensolution. It does not imply evaluation before solving the reduced problems; residual-based a posteriori estimators can also avoid the full-order eigensolution \cite{jakobsson2011posteriori}. The CB decomposition below establishes the notation and the connection to the rank-1 HCB-1 formulation.

\subsection{Error measures and definitions}\label{sec:estimator}
Let $(\tilde{\lambda}_i, \tilde{\boldsymbol{\phi}}_i)$ denote the $i$-th eigenpair of the CB reduced eigenvalue problem
\begin{align}\label{eq:cb_reduced_eig}
    \widetilde{\mathbf{K}}_{\mathrm{CB}} \tilde{\boldsymbol{\phi}}_i = \tilde{\lambda}_i \widetilde{\mathbf{M}}_{\mathrm{CB}} \tilde{\boldsymbol{\phi}}_i,
\end{align}
mass-normalized such that $\tilde{\boldsymbol{\phi}}_i^{\mathrm{T}} \widetilde{\mathbf{M}}_{\mathrm{CB}} \tilde{\boldsymbol{\phi}}_j = \delta_{ij}$ and $\tilde{\boldsymbol{\phi}}_i^{\mathrm{T}} \widetilde{\mathbf{K}}_{\mathrm{CB}} \tilde{\boldsymbol{\phi}}_j = \tilde{\lambda}_i \delta_{ij}$. The standard CB approximation of the $i$-th global eigenvector is $\mathbf{T}_{\mathrm{CB}} \tilde{\boldsymbol{\phi}}_i$, which neglects the contribution of the truncated modes.

Let $(\lambda_i, \boldsymbol{\phi}_i)$ denote the $i$-th exact eigenpair from Eq.~\eqref{eq:eigproblem}. The relative eigenvalue error of the $i$-th CB eigenvalue is defined as
\begin{align}\label{eq:rel_error_def}
    \xi_i \coloneqq \frac{\tilde{\lambda}_i - \lambda_i}{\lambda_i},
\end{align}
where $\lambda_i$ is the exact eigenvalue. The auxiliary matrix used throughout the error estimator derivations is
\begin{align}\label{eq:Gi}
    \mathbf{G}(\mu) \coloneqq \mathbf{M}_p - \frac{1}{\mu} \mathbf{K}_p.
\end{align}

Suppose an enhanced method produces eigenpairs $(\bar{\lambda}_i, \bar{\boldsymbol{\phi}}_i)$ in the reduced space. The corresponding physical eigenvector is decomposed as
\begin{align}\label{eq:decomp_general}
    \mathbf{u}_i = (\mathbf{T}_0 + \mathbf{T}_r) \bar{\boldsymbol{\phi}}_i = \mathbf{u}_{0,i} + \mathbf{u}_{r,i},
\end{align}
where $\mathbf{T}_0$ represents the base (CB-like) contribution and $\mathbf{T}_r$ represents the residual correction. A Rayleigh quotient perturbation analysis of the base-only approximation (Appendix~\ref{sec:derivation_decomp}) gives the following leading-order expression for the supporting CB error estimate:
\begin{align}\label{eq:eta_enhanced}
    \hat{\eta}_i = 2\bar{\boldsymbol{\phi}}_i^{\mathrm{T}} \mathbf{T}_0^{\mathrm{T}} \mathbf{G}(\bar{\lambda}_i) \, \mathbf{T}_r \bar{\boldsymbol{\phi}}_i + \bar{\boldsymbol{\phi}}_i^{\mathrm{T}} \mathbf{T}_r^{\mathrm{T}} \mathbf{G}(\bar{\lambda}_i) \, \mathbf{T}_r \bar{\boldsymbol{\phi}}_i,
\end{align}
where $\hat{\eta}_i$ is the estimated relative eigenvalue error, $\bar{\lambda}_i$ is the $i$-th eigenvalue of the enhanced (HCB) reduced system, and $\mathbf{G}(\bar{\lambda}_i) = \mathbf{M}_p - \bar{\lambda}_i^{-1} \mathbf{K}_p$ is the auxiliary matrix (Eq.~\eqref{eq:Gi}). Its interpretation as an estimate of $\xi_i$ additionally requires $\mathbf{u}_{0,i}^{\mathrm{T}}\mathbf{M}_p\mathbf{u}_{0,i}\approx 1$ and a base-only Rayleigh quotient close to the corresponding CB eigenvalue. Both the reference eigenvalue error and the difference between the base-only Rayleigh quotient and the CB eigenvalue must be small on the error scale being estimated.

\subsection{Extension of the CB error estimator via HCB decomposition}\label{sec:estimator_enhanced}\label{sec:estimator_hcb_detail}

As a supporting extension of the CB estimator in \cite{kim2014estimating}, the physical eigenvector of the HCB-$n$ method $\mathbf{u}_i = \mathbf{T}_{\mathrm{HCB\text{-}}n} \mathbf{T}_{\mathrm{SEREP}} \bar{\boldsymbol{\phi}}_i$ is decomposed into a CB-like component and a residual correction component following Eq.~\eqref{eq:decomp_general}. The composite transformation is decomposed by partitioning the SEREP matrix according to the CB and residual DOFs:
\begin{align}\label{eq:serep_partition}
    \mathbf{T}_{\mathrm{SEREP}} = \begin{bmatrix} \mathbf{V}_{\mathrm{CB}} \\ \mathbf{V}_{\mathrm{res}} \end{bmatrix},
\end{align}
where $\mathbf{V}_{\mathrm{CB}} \in \mathbb{R}^{n_{\mathrm{CB}} \times n_{\mathrm{CB}}}$ corresponds to the CB generalized coordinates ($\boldsymbol{\Phi}_d^{(k)}$, $\boldsymbol{\Psi}_b^{(k)}$, $\mathbf{u}_b$) and $\mathbf{V}_{\mathrm{res}}$ corresponds to the residual mode coordinates ($\boldsymbol{\eta}_1$ for HCB-1; $\boldsymbol{\eta}_1$, $\boldsymbol{\eta}_2$ for HCB-2). Since the first $n_{\mathrm{CB}}$ columns of $\mathbf{T}_{\mathrm{HCB\text{-}}n}$ are $\mathbf{T}_{\mathrm{CB}}$ and the remaining columns contain the residual modes, the base and residual components are
\begin{align}\label{eq:T0_Tr_hcb}
    \mathbf{T}_0 \coloneqq \mathbf{T}_{\mathrm{CB}} \mathbf{V}_{\mathrm{CB}}, \qquad \mathbf{T}_r \coloneqq [\mathbf{R}_1, \ldots, \mathbf{R}_n] \mathbf{V}_{\mathrm{res}},
\end{align}
where $[\mathbf{R}_1, \ldots, \mathbf{R}_n]$ denotes the assembled residual mode columns of $\mathbf{T}_{\mathrm{HCB\text{-}}n}$. The factor $\mathbf{V}_{\mathrm{CB}}$ maps retained SEREP modal coordinates to CB coordinates and is generally not the identity matrix, even when the residual contribution is small. The supporting CB error estimator is
\begin{align}\label{eq:eta_hcb_n}
    \hat{\eta}_{\mathrm{CB},i}^{(n)} = 2\bar{\boldsymbol{\phi}}_i^{\mathrm{T}} \mathbf{T}_0^{\mathrm{T}} \mathbf{G}(\bar{\lambda}_i) \, \mathbf{T}_r \bar{\boldsymbol{\phi}}_i + \bar{\boldsymbol{\phi}}_i^{\mathrm{T}} \mathbf{T}_r^{\mathrm{T}} \mathbf{G}(\bar{\lambda}_i) \, \mathbf{T}_r \bar{\boldsymbol{\phi}}_i,
\end{align}
where $(\bar{\lambda}_i, \bar{\boldsymbol{\phi}}_i)$ is the $i$-th HCB-$n$ eigenpair after SEREP reduction. For computational efficiency, the cross-product matrices
\begin{align}
    \mathbf{M}_{0r} = \mathbf{T}_0^{\mathrm{T}} \mathbf{M}_p \mathbf{T}_r, \quad \mathbf{K}_{0r} = \mathbf{T}_0^{\mathrm{T}} \mathbf{K}_p \mathbf{T}_r, \quad \mathbf{M}_{rr} = \mathbf{T}_r^{\mathrm{T}} \mathbf{M}_p \mathbf{T}_r, \quad \mathbf{K}_{rr} = \mathbf{T}_r^{\mathrm{T}} \mathbf{K}_p \mathbf{T}_r
\end{align}
are precomputed once, and the estimator for each mode is evaluated as
\begin{empheq}[box=\fbox]{align}\label{eq:eta_hcb_efficient}
    \hat{\eta}_{\mathrm{CB},i}^{(n)} = 2\bar{\boldsymbol{\phi}}_i^{\mathrm{T}} \left(\mathbf{M}_{0r} - \frac{1}{\bar{\lambda}_i} \mathbf{K}_{0r}\right) \bar{\boldsymbol{\phi}}_i + \bar{\boldsymbol{\phi}}_i^{\mathrm{T}} \left(\mathbf{M}_{rr} - \frac{1}{\bar{\lambda}_i} \mathbf{K}_{rr}\right) \bar{\boldsymbol{\phi}}_i.
\end{empheq}

Both $\mathbf{T}_0$ and $\mathbf{T}_r$ depend on the HCB order through the corresponding SEREP partition. Thus, each estimator $\hat{\eta}_{\mathrm{CB},i}^{(n)}$ must use the partition and eigenpair of the selected HCB-$n$ model. Improved reference eigenvalues alone do not guarantee improved CB error estimates, because the base-only Rayleigh quotient also changes.

\begin{remark}[Exact-denominator variant]
The derivation of Eq.~\eqref{eq:eta_enhanced} in Appendix~\ref{sec:derivation_decomp} invokes the leading-order approximation $\mathbf{u}_{0,i}^{\mathrm{T}} \mathbf{M}_p \mathbf{u}_{0,i} \approx 1$, which holds when the residual component $\mathbf{u}_{r,i} = \mathbf{T}_r \bar{\boldsymbol{\phi}}_i$ is small relative to the base component. The approximation can be eliminated entirely by introducing the additional precomputed matrices
\begin{align*}
    \mathbf{M}_{00} = \mathbf{T}_0^{\mathrm{T}} \mathbf{M}_p \mathbf{T}_0, \qquad \mathbf{K}_{00} = \mathbf{T}_0^{\mathrm{T}} \mathbf{K}_p \mathbf{T}_0,
\end{align*}
which are of the same computational class as $\mathbf{M}_{0r}$ and $\mathbf{K}_{0r}$ in Eq.~\eqref{eq:eta_hcb_efficient} (additional matrix products of the same form). The exact-denominator estimator is then
\begin{align*}
    \hat{\eta}_{\mathrm{CB},i}^{\mathrm{ex}} = \frac{\bar{\boldsymbol{\phi}}_i^{\mathrm{T}} \mathbf{K}_{00} \bar{\boldsymbol{\phi}}_i - \bar{\lambda}_i \, \bar{\boldsymbol{\phi}}_i^{\mathrm{T}} \mathbf{M}_{00} \bar{\boldsymbol{\phi}}_i}{\bar{\lambda}_i \, \bar{\boldsymbol{\phi}}_i^{\mathrm{T}} \mathbf{M}_{00} \bar{\boldsymbol{\phi}}_i},
\end{align*}
which evaluates the Rayleigh quotient of the base-only approximation $\mathbf{u}_{0,i}$ relative to $\bar{\lambda}_i$ without any leading-order approximation. The relation $\hat{\eta}_{\mathrm{CB},i}^{\mathrm{ex}} \approx \hat{\eta}_{\mathrm{CB},i}^{(n)}$ holds when $\bar{\boldsymbol{\phi}}_i^{\mathrm{T}} \mathbf{M}_{00} \bar{\boldsymbol{\phi}}_i \approx 1$, confirming that Eq.~\eqref{eq:eta_hcb_efficient} is the leading-order approximation of $\hat{\eta}_{\mathrm{CB},i}^{\mathrm{ex}}$. This identity removes the unit-mass approximation for the base component; it does not by itself establish the accuracy of the resulting CB error estimate.
\end{remark}

\subsection{HCB-1 error estimator via HCB-2}\label{sec:estimator_hcb1}

Before SEREP reduction, the augmented transformation matrices form a nested hierarchy:\label{sec:morth}
\begin{align}\label{eq:nested}
    \mathcal{R}(\mathbf{T}_{\mathrm{CB}}) \subseteq \mathcal{R}(\mathbf{T}_{\mathrm{HCB\text{-}1}}) \subseteq \mathcal{R}(\mathbf{T}_{\mathrm{HCB\text{-}2}}) \subseteq \mathbb{R}^{N_g},
\end{align}
where $\mathcal{R}(\cdot)$ denotes the column space. The spaces may coincide if newly added vectors are linearly dependent. The following argument assumes that each space is represented by a full-column-rank basis and that the resulting Ritz problems are solved in exact arithmetic. These are mathematical assumptions; column normalization alone does not establish linear independence or numerical stability. For symmetric $\mathbf{K}_p$ and positive-definite $\mathbf{M}_p$, the Courant--Fischer min-max principle \cite{courant1953methods, parlett1998symmetric, golub2013matrix} orders the Ritz eigenvalues by their common index in ascending order. SEREP retains the selected augmented eigenvectors and preserves their eigenvalues in exact arithmetic. Thus, for the retained eigenvalues,
\begin{align}\label{eq:minmax}
    \lambda_{\mathrm{CB},i} \geq \lambda_{\mathrm{HCB},i}^{(1)} \geq \lambda_{\mathrm{HCB},i}^{(2)} \geq \lambda_i,
\end{align}
where $\lambda_i$ denotes the exact eigenvalue with the same ordered index. The relative errors below are restricted to $\lambda_i>0$. This ordering guarantees non-negative differences. If mode correlation is used to pair closely spaced modes, the pairing must also respect the ordered indices for this lower-bound argument to apply. The bound is a statement about ordered Ritz eigenvalues; it does not depend on how well the individual eigenvectors align.

To apply the general decomposition framework (Eq.~\eqref{eq:eta_enhanced}) for estimating the HCB-1 error, each HCB-2 physical eigenvector must be decomposed into an HCB-1-like component and a residual. Let $(\lambda_{\mathrm{HCB},i}^{(1)}, \bar{\boldsymbol{\phi}}_i^{(1)})$ and $(\lambda_{\mathrm{HCB},i}^{(2)}, \bar{\boldsymbol{\phi}}_i^{(2)})$ denote the $i$-th eigenpairs of the HCB-1 and HCB-2 reduced systems after SEREP reduction, respectively. Define the mode correspondence coefficient
\begin{align}\label{eq:mode_corr}
    c_i \coloneqq (\bar{\boldsymbol{\phi}}_i^{(1)})^{\mathrm{T}} \mathbf{T}_{\mathrm{comp\text{-}1}}^{\mathrm{T}} \mathbf{M}_p \, \mathbf{T}_{\mathrm{comp\text{-}2}} \bar{\boldsymbol{\phi}}_i^{(2)},
\end{align}
where $\mathbf{T}_{\mathrm{comp\text{-}}n} = \mathbf{T}_{\mathrm{HCB\text{-}}n} \mathbf{T}_{\mathrm{SEREP\text{-}}n}$ is the composite transformation that maps reduced coordinates to physical coordinates for the HCB-$n$ method. The coefficient $c_i$ measures the $\mathbf{M}_p$-inner-product alignment between the $i$-th HCB-1 and HCB-2 eigenmodes. For mass-normalized physical eigenvectors, $0\leq c_i^2\leq 1$. Well-aligned modes have $|c_i|\approx 1$; the sign of $c_i$ depends on the arbitrary eigenvector signs.

The rank-1 mass-orthogonal projection of the $i$-th HCB-2 physical eigenvector onto the $i$-th HCB-1 mode defines the decomposition\label{sec:direct_diff}
\begin{align}\label{eq:rank1_decomp}
    \mathbf{u}_{0,i} \coloneqq c_i \, \mathbf{T}_{\mathrm{comp\text{-}1}} \bar{\boldsymbol{\phi}}_i^{(1)}, \qquad
    \mathbf{u}_{r,i} \coloneqq \mathbf{T}_{\mathrm{comp\text{-}2}} \bar{\boldsymbol{\phi}}_i^{(2)} - \mathbf{u}_{0,i}.
\end{align}
The above decomposition makes the base and residual components mass-orthogonal:
\begin{align}
    \mathbf{u}_{0,i}^{\mathrm{T}} \mathbf{M}_p \, \mathbf{u}_{r,i} = c_i\bigl(c_i - c_i\bigr) = 0,
\end{align}
using the $\mathbf{M}_p$-normalization of the physical HCB eigenvectors. The base component has the quadratic forms
\begin{align}\label{eq:rank1_quad}
    \mathbf{u}_{0,i}^{\mathrm{T}} \mathbf{M}_p \, \mathbf{u}_{0,i} = c_i^2, \qquad
    \mathbf{u}_{0,i}^{\mathrm{T}} \mathbf{K}_p \, \mathbf{u}_{0,i} = c_i^2 \lambda_{\mathrm{HCB},i}^{(1)},
\end{align}
since the physical vector $\mathbf{T}_{\mathrm{comp\text{-}1}}\bar{\boldsymbol{\phi}}_i^{(1)}$ is mass-normalized and has Rayleigh quotient $\lambda_{\mathrm{HCB},i}^{(1)}$.

Substituting into the general error estimator (Eq.~\eqref{eq:eta_enhanced}) and expanding using $\mathbf{u}_i^{\mathrm{T}} \mathbf{K}_p \mathbf{u}_i = \bar{\lambda}_i$, $\mathbf{u}_i^{\mathrm{T}} \mathbf{M}_p \mathbf{u}_i = 1$, and Eq.~\eqref{eq:rank1_quad}, all terms involving the stiffness cross-product $\mathbf{u}_{0,i}^{\mathrm{T}} \mathbf{K}_p \mathbf{u}_{r,i}$ cancel exactly, yielding
\begin{align}\label{eq:eta_hcb1_morth}
    \hat{\eta}_{\mathrm{proj},i}^{(1)} = c_i^2 \, \frac{\lambda_{\mathrm{HCB},i}^{(1)} - \bar{\lambda}_i}{\bar{\lambda}_i}.
\end{align}
This is an algebraic identity for the projection-based, leading-order expression in Eq.~\eqref{eq:eta_enhanced}, not an identity for the full-order eigenvalue error. With $\bar{\lambda}_i=\lambda_{\mathrm{HCB},i}^{(2)}$, define the direct eigenvalue-difference estimator separately as
\begin{empheq}[box=\fbox]{align}\label{eq:eta_direct}
    \hat{\eta}_{\mathrm{HCB},i}^{(1)} \coloneqq \frac{\lambda_{\mathrm{HCB},i}^{(1)} - \lambda_{\mathrm{HCB},i}^{(2)}}{\lambda_{\mathrm{HCB},i}^{(2)}}.
\end{empheq}
The relation $\hat{\eta}_{\mathrm{proj},i}^{(1)}=c_i^2\hat{\eta}_{\mathrm{HCB},i}^{(1)}$ is exact. Replacing the projection expression by the direct estimator uses $c_i^2\approx 1$; this approximation only concerns the agreement of these two expressions. It is not used to define the direct estimator or to prove its lower-bound property. All HCB-1 estimates reported in Section~\ref{sec:example} use Eq.~\eqref{eq:eta_direct} without the factor $c_i^2$. Writing $\lambda_{\mathrm{HCB},i}^{(1)} = \lambda_i(1 + \xi_{\mathrm{HCB},i}^{(1)})$ and $\lambda_{\mathrm{HCB},i}^{(2)} = \lambda_i(1 + \xi_{\mathrm{HCB},i}^{(2)})$, the estimator satisfies the exact identity
\begin{align}\label{eq:direct_decomposition}
    \hat{\eta}_{\mathrm{HCB},i}^{(1)}
    = \frac{\lambda_{\mathrm{HCB},i}^{(1)} - \lambda_{\mathrm{HCB},i}^{(2)}}{\lambda_{\mathrm{HCB},i}^{(2)}}
    = \frac{\xi_{\mathrm{HCB},i}^{(1)} - \xi_{\mathrm{HCB},i}^{(2)}}{1 + \xi_{\mathrm{HCB},i}^{(2)}},
\end{align}
or equivalently,
\begin{align}\label{eq:direct_approx}
    \xi_{\mathrm{HCB},i}^{(1)} = \hat{\eta}_{\mathrm{HCB},i}^{(1)} + \xi_{\mathrm{HCB},i}^{(2)} + \hat{\eta}_{\mathrm{HCB},i}^{(1)} \, \xi_{\mathrm{HCB},i}^{(2)}.
\end{align}
Under Eq.~\eqref{eq:minmax} and $\lambda_i>0$, Eq.~\eqref{eq:direct_decomposition} gives
\begin{align*}
    0\leq \hat{\eta}_{\mathrm{proj},i}^{(1)}
    \leq \hat{\eta}_{\mathrm{HCB},i}^{(1)}
    \leq \xi_{\mathrm{HCB},i}^{(1)}.
\end{align*}
For $\xi_{\mathrm{HCB},i}^{(1)}>0$, the estimation-to-error ratio is
\begin{align*}
    \frac{\hat{\eta}_{\mathrm{HCB},i}^{(1)}}{\xi_{\mathrm{HCB},i}^{(1)}}
    = \frac{1-\xi_{\mathrm{HCB},i}^{(2)}/\xi_{\mathrm{HCB},i}^{(1)}}{1+\xi_{\mathrm{HCB},i}^{(2)}}.
\end{align*}
A close estimate therefore requires both $\xi_{\mathrm{HCB},i}^{(2)}\ll 1$ and $\xi_{\mathrm{HCB},i}^{(2)}/\xi_{\mathrm{HCB},i}^{(1)}\ll 1$. Neither monotone ordering nor $c_i^2\approx 1$ guarantees these conditions. A small lower bound consequently does not certify that the actual error is below a prescribed tolerance. The conditions characterize estimator accuracy in terms of the unknown full-order errors; the present formulation does not provide a computable upper bound or an independent test that certifies them without a reference solution. The direct estimator requires only the HCB-1 and HCB-2 eigenvalues; evaluating $c_i$ is needed to assess mode alignment and the projection-based expression, but is not part of Eq.~\eqref{eq:eta_direct}.

\begin{remark}[Properties of the mass-orthogonal projection estimator]\label{rmk:morth}
    For $c_i\neq 0$, the rank-1 projection makes the Rayleigh quotient of $\mathbf{u}_{0,i}$ equal to $\lambda_{\mathrm{HCB},i}^{(1)}$. Projection onto the entire retained HCB-1 space generally mixes several modes, whose combined Rayleigh quotient need not equal the $i$-th eigenvalue. Rank-1 projection therefore provides a mode-specific connection between Eq.~\eqref{eq:eta_enhanced} and the direct eigenvalue difference. Closely spaced or repeated eigenvalues require care in interpreting individual mode correspondence.
\end{remark}

\section{Numerical examples}\label{sec:example}

Three numerical examples of increasing geometric complexity are presented to validate the proposed error estimators. All examples use the same material properties: Young's modulus $E = 210$~GPa, Poisson's ratio $\nu = 0.3$, and density $\rho = 7{,}850$~kg/m$^3$. For each example, three CMS methods---CB, HCB-1, and HCB-2---are applied at the same reduced-model size, and the first 20 eigenvalues are compared against the full-order model (FOM) solution.

The two estimators are evaluated as follows:
\begin{enumerate}
    \item HCB-1 error estimator (HCB-2-based): The direct estimator $\hat{\eta}_{\mathrm{HCB},i}^{(1)}$ in Eq.~\eqref{eq:eta_direct} estimates the HCB-1 relative eigenvalue error $\xi_{\mathrm{HCB},i}^{(1)}$. The reported values exclude the factor $c_i^2$ in the projection-based expression.
    \item Supporting CB error estimator (HCB-1-based): The estimator $\hat{\eta}_{\mathrm{CB},i}^{(1)}$ in Eq.~\eqref{eq:eta_hcb_efficient} with $n=1$ uses HCB-1 eigenpairs and the SEREP partition to estimate the CB relative eigenvalue error $\xi_i$.
\end{enumerate}
The estimation quality is quantified by the ratio $\hat{\eta} / \xi$; a ratio of 1.0 indicates a perfect estimate.

\begin{remark}[Mode correspondence and MAC]\label{sec:error}
    For real mode shapes, an unsquared form of the modal assurance criterion (MAC) \cite{pastor2012modal} is used in this study:
    \begin{align*}
        \mathrm{MAC}_{ij} &=
        \frac{|(\boldsymbol{\phi}_i^{\mathrm{exact}})^{\mathrm{T}}\boldsymbol{\phi}_j^{\mathrm{phys}}|}
        {\|\boldsymbol{\phi}_i^{\mathrm{exact}}\|_2\,\|\boldsymbol{\phi}_j^{\mathrm{phys}}\|_2}.
    \end{align*}
    Here, $\boldsymbol{\phi}_j^{\mathrm{phys}}$ is the CB or HCB eigenvector expanded to physical coordinates, and $\|\cdot\|_2$ denotes the Euclidean norm. Values range from 0 for orthogonal modes to 1 for aligned modes. Mode correspondence is identified using the largest MAC values. This Euclidean correlation differs from the mass-weighted coefficient $c_i$ in Eq.~\eqref{eq:mode_corr}.
\end{remark}

\subsection{Example 1: Cantilever plate}\label{sec:ex_plate}
\begin{figure}[ht!]
\centering
\includegraphics[width=0.85\textwidth]{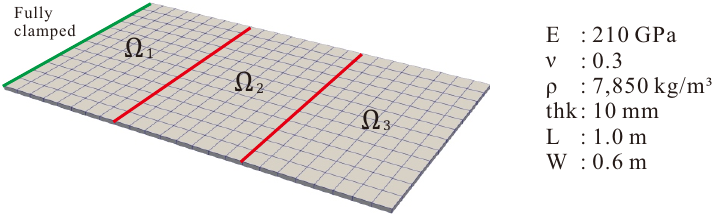}
\caption{Cantilever plate example. Geometry with the clamped boundary at $x = 0$ and three-substructure partition along the $x$-direction at $x = L/3$ and $x = 2L/3$. The plate ($1.0 \times 0.6 \times 0.01$~m) is meshed with a uniform $24\times12\times1$ grid of 8-node hexahedral (hex8) elements (650 nodes, 1{,}872 free DOFs). Substructures~1 and~2 each occupy one-third of the plate length, while Substructure~3 covers the remaining third.}
\label{fig:plate_mesh}
\end{figure}

A cantilever rectangular plate ($1.0 \times 0.6 \times 0.01$~m) is discretized with a $24 \times 12 \times 1$ mesh of 8-node hexahedral elements (650 nodes, 1{,}872 free DOFs after clamping the $x = 0$ face). The structure is partitioned into three substructures along the $x$-direction at $x = L/3$ and $x = 2L/3$, as summarized in Table~\ref{tab:plate_subs} and illustrated in Fig.~\ref{fig:plate_mesh}. The CB reduced-model size is $n_{\mathrm{CB}} = 10 + 10 + 8 + 156 = 184$ (9.8\% of 1{,}872 DOFs).

\begin{table}[!htbp]
\centering
\caption{Substructure configuration for the cantilever plate.}
\label{tab:plate_subs}
\begin{tabular}{cccc}
\hline
Substructure & Interior DOFs ($N_s^{(k)}$) & Boundary DOFs ($N_b$) & Retained modes ($N_d^{(k)}$) \\
\hline
1 & 546 & -- & 10 \\
2 & 546 & -- & 10 \\
3 & 624 & -- & 8 \\
Boundary & -- & 156 & -- \\
\hline
\end{tabular}
\end{table}

\begin{figure}[ht!]
\centering
\includegraphics[width=\textwidth]{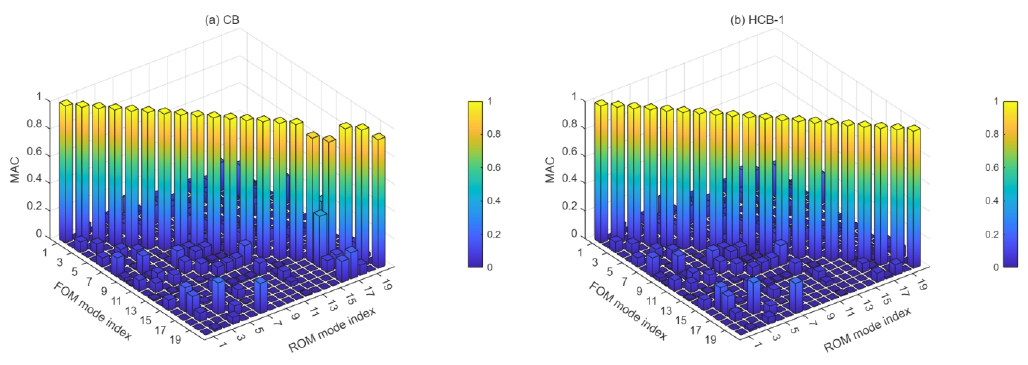}
\caption{Modal assurance criterion $\mathrm{MAC}_{ij}$ for the cantilever plate. (a)~CB and (b)~HCB-1 modes compared with FOM modes for the first 20 modes. MAC is defined in Remark~\ref{sec:error}; values of 1 and 0 indicate aligned and orthogonal modes, respectively.}
\label{fig:plate_mac}
\end{figure}

\begin{figure}[ht!]
\centering
\includegraphics[width=\textwidth]{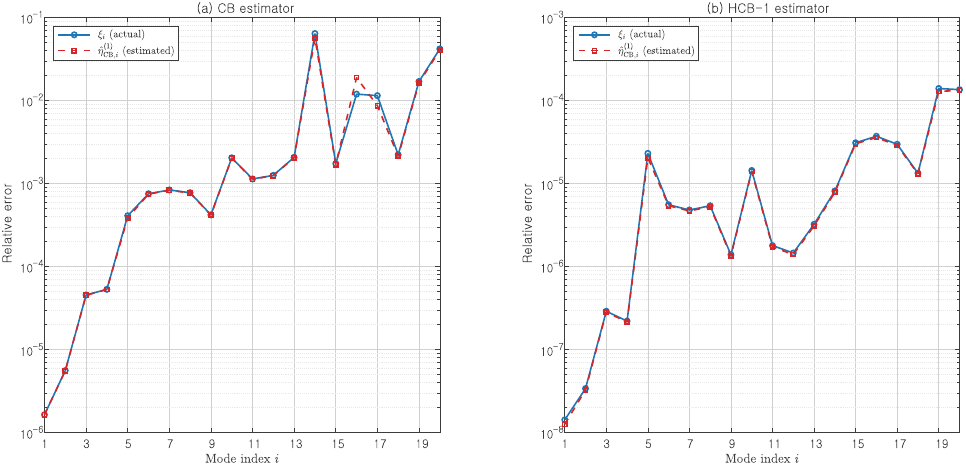}
\caption{Error estimator validation for the cantilever plate. Actual relative eigenvalue error $\xi_i$ (solid blue circles) and estimated error $\hat{\eta}_i$ (dashed red squares) plotted on a logarithmic scale for all 20 modes. (a)~CB estimator $\hat{\eta}_{\mathrm{CB},i}^{(1)}$, with ratios in the range $0.76$--$1.58$ across all modes. (b)~HCB-1 estimator $\hat{\eta}_{\mathrm{HCB},i}^{(1)}$, showing a consistent underestimate with ratios $0.88$--$0.99$. The HCB-1 errors are smaller than the CB errors for all 20 modes.}
\label{fig:plate_err}
\end{figure}

The MAC matrices in Fig.~\ref{fig:plate_mac} show the correspondence between FOM and reduced modes, as defined in Remark~\ref{sec:error}. The actual and estimated errors for all modes are compared in Fig.~\ref{fig:plate_err}. For the first 10 modes listed in Table~\ref{tab:plate_results}, the CB error estimator achieves ratios between 0.94 and 1.00, indicating close agreement with the actual CB errors ($\sim 10^{-6}$--$10^{-3}$). The HCB-1 estimator achieves ratios between 0.88 and 0.98, providing a consistent slight underestimate of the HCB-1 errors ($\sim 10^{-8}$--$10^{-5}$). The HCB-1 errors are smaller than the CB errors for all 20 modes, demonstrating the accuracy improvement from the residual mode correction. The estimation ratio of approximately 0.885, observed for modes 1 and 5, admits a precise interpretation through Eq.~\eqref{eq:direct_decomposition}: rearranging gives
\begin{align}
    \frac{\hat{\eta}_{\mathrm{HCB},i}^{(1)}}{\xi_{\mathrm{HCB},i}^{(1)}}
    = \frac{1}{1 + \xi_{\mathrm{HCB},i}^{(2)}/\hat{\eta}_{\mathrm{HCB},i}^{(1)} + \xi_{\mathrm{HCB},i}^{(2)}},
\end{align}
For these small-error cases, the final term in the denominator is negligible, so the estimation ratio degrades when $\xi_{\mathrm{HCB},i}^{(2)}$ is non-negligible relative to $\hat{\eta}_{\mathrm{HCB},i}^{(1)}$. For modes 1 and 5, $\xi_{\mathrm{HCB\text{-}2}}/\hat{\eta} \approx 0.13$, whereas for the remaining modes reported in Table~\ref{tab:plate_results} it is $0.02$--$0.04$; accordingly, the ratio is $1/(1+0.13) \approx 0.885$ for modes 1 and 5 and $0.96$--$0.98$ for the others. This larger $\xi_{\mathrm{HCB\text{-}2}}/\hat{\eta}$ ratio is equivalent to a smaller improvement factor $\xi_{\mathrm{HCB\text{-}1}}/\xi_{\mathrm{HCB\text{-}2}}$ from the first to the second order of the HCB series: for modes 1 and 5, $\xi_{\mathrm{HCB\text{-}1}}/\xi_{\mathrm{HCB\text{-}2}} \approx 8.7$, compared to $26$--$42$ for modes 2--4 and 6--10. The slower per-step convergence indicates that the second-order residual correction is relatively less effective at reducing the HCB-1 error for these two modes. Comparison with the other examples supports this interpretation: for the elbow pipe, modes 6 and 7 exhibit $\xi_{\mathrm{HCB\text{-}1}}/\xi_{\mathrm{HCB\text{-}2}} \approx 21$--$25$, yielding estimation ratios of $0.953$--$0.961$; for the RPV, the minimum per-step factor is $\approx 6.6$ (mode 3) and ratios range from $0.849$ to $0.989$. These observations are consistent with the error relation in Eq.~\eqref{eq:direct_decomposition}: close estimation requires the HCB-2 error to be small both in magnitude and relative to the HCB-1 error.

\begin{table}[!htbp]
\centering
\caption{Error estimator validation for the cantilever plate (first 10 of 20 modes).}
\label{tab:plate_results}
\begin{tabular}{cccccccc}
\hline
\multicolumn{2}{c}{} & \multicolumn{3}{c}{CB} & \multicolumn{3}{c}{HCB-1} \\
\cline{3-5}\cline{6-8}
Mode & $f_{\mathrm{FOM}}$ [Hz] & $\xi$ & $\hat{\eta}$ & ratio & $\xi$ & $\hat{\eta}$ & ratio \\
\hline
1  &   23.64 & 1.64e-06 & 1.65e-06 & 1.000 & 1.43e-08 & 1.26e-08 & 0.885 \\
2  &   41.65 & 5.59e-06 & 5.56e-06 & 0.994 & 3.40e-08 & 3.27e-08 & 0.962 \\
3  &  148.31 & 4.55e-05 & 4.52e-05 & 0.993 & 2.90e-07 & 2.83e-07 & 0.976 \\
4  &  176.38 & 5.33e-05 & 5.31e-05 & 0.996 & 2.21e-07 & 2.16e-07 & 0.976 \\
5  &  409.46 & 4.10e-04 & 3.87e-04 & 0.943 & 2.31e-05 & 2.05e-05 & 0.885 \\
6  &  416.63 & 7.54e-04 & 7.47e-04 & 0.991 & 5.56e-06 & 5.42e-06 & 0.974 \\
7  &  442.62 & 8.39e-04 & 8.33e-04 & 0.992 & 4.81e-06 & 4.69e-06 & 0.974 \\
8  &  496.75 & 7.77e-04 & 7.71e-04 & 0.993 & 5.41e-06 & 5.26e-06 & 0.972 \\
9  &  549.93 & 4.21e-04 & 4.18e-04 & 0.994 & 1.39e-06 & 1.35e-06 & 0.973 \\
10 &  718.14 & 2.05e-03 & 2.02e-03 & 0.988 & 1.44e-05 & 1.40e-05 & 0.971 \\
\hline
\end{tabular}
\end{table}

\subsection{Example 2: 90-degree elbow pipe}\label{sec:ex_pipe}
\begin{figure}[ht!]
\centering
\includegraphics[width=0.85\textwidth]{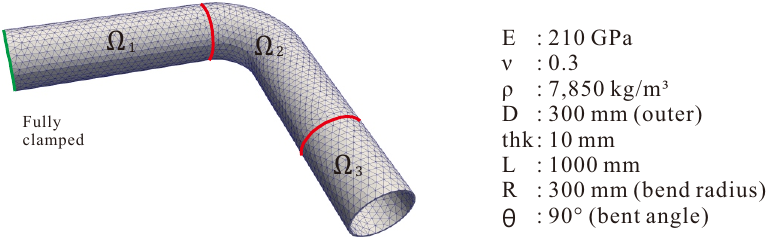}
\caption{90-degree elbow pipe example. Geometry with clamped ends and three-substructure partition. Substructure~1 ($\Omega_1$) covers the horizontal straight section ($L = 1.0$~m), Substructure~2 ($\Omega_2$) covers the 90-degree bend and the upper vertical section, and Substructure~3 ($\Omega_3$) covers the lower vertical section ($L/2 = 0.5$~m). The pipe has outer radius $R_o = 0.15$~m, inner radius $R_i = 0.14$~m, and bend centerline radius $R_{\mathrm{bend}} = 0.3$~m; the unstructured tet4 mesh has 14{,}092 elements (4{,}653 nodes, 13{,}641 free DOFs).}
\label{fig:pipe_mesh}
\end{figure}

A 90-degree elbow pipe consisting of a horizontal straight section ($L = 1.0$~m), a 90-degree bend (centerline radius $R_{\mathrm{bend}} = 0.3$~m), and a vertical straight section ($L = 1.0$~m) is considered. The pipe cross-section is a hollow annulus with outer radius $R_o = 0.15$~m and inner radius $R_i = 0.14$~m (wall thickness $t = 0.01$~m). Both ends are clamped. The geometry is meshed with 14{,}092 four-node tetrahedral (tet4) elements (4{,}653 nodes, 13{,}641 free DOFs), as illustrated in Fig.~\ref{fig:pipe_mesh}.

The structure is partitioned into three substructures by element centroid coordinates: Substructure~1 ($\Omega_1$) covers the horizontal straight section ($x < L$), Substructure~2 ($\Omega_2$) covers the 90-degree bend and the upper vertical section, and Substructure~3 ($\Omega_3$) covers the lower vertical section. Table~\ref{tab:pipe_subs} summarizes the configuration. The CB reduced-model size is $n_{\mathrm{CB}} = 15 + 15 + 15 + 321 = 366$ (2.7\% of 13{,}641 DOFs).

\begin{table}[!htbp]
\centering
\caption{Substructure configuration for the 90-degree elbow pipe.}
\label{tab:pipe_subs}
\begin{tabular}{cccc}
\hline
Substructure & Interior DOFs ($N_s^{(k)}$) & Boundary DOFs ($N_b$) & Retained modes ($N_d^{(k)}$) \\
\hline
1 (horizontal straight) & 5{,}325 & -- & 15 \\
2 (arc + upper vertical) & 5{,}340 & -- & 15 \\
3 (lower vertical) & 2{,}655 & -- & 15 \\
Boundary & -- & 321 & -- \\
\hline
\end{tabular}
\end{table}

\begin{figure}[ht!]
\centering
\includegraphics[width=\textwidth]{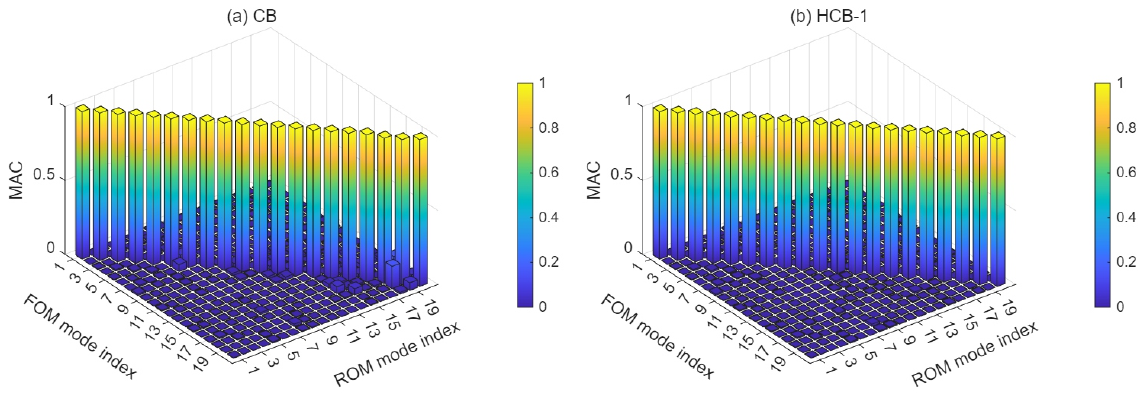}
\caption{Modal assurance criterion $\mathrm{MAC}_{ij}$ for the 90-degree elbow pipe. (a)~CB and (b)~HCB-1 modes compared with FOM modes for the first 20 modes. MAC is defined in Remark~\ref{sec:error}; values of 1 and 0 indicate aligned and orthogonal modes, respectively.}
\label{fig:pipe_mac}
\end{figure}

\begin{figure}[ht!]
\centering
\includegraphics[width=\textwidth]{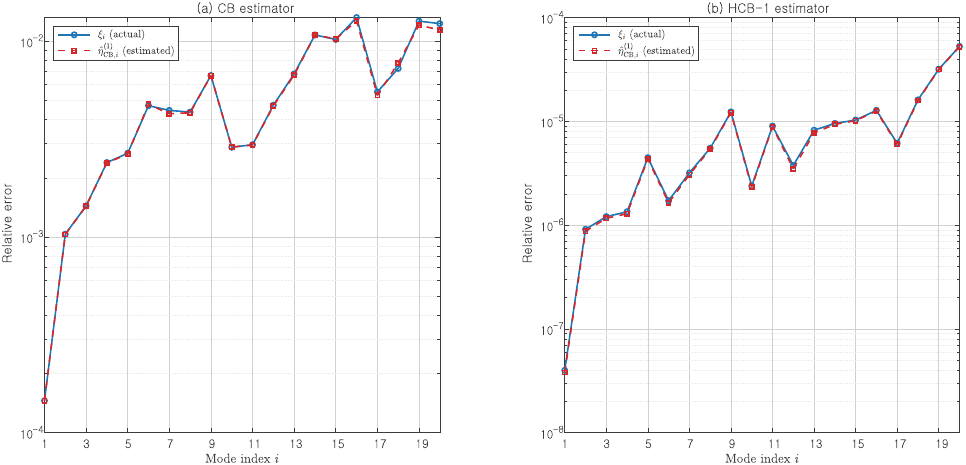}
\caption{Error estimator validation for the 90-degree elbow pipe. Actual relative eigenvalue error $\xi_i$ (solid blue circles) and estimated error $\hat{\eta}_i$ (dashed red squares) for all 20 modes. (a)~CB estimator, with ratios in the range $0.93$--$1.06$ and slight overestimation observed for mode~18 (ratio~$\approx 1.06$). (b)~HCB-1 estimator, with ratios $0.93$--$0.99$. The HCB-1 errors are smaller than the CB errors for all 20 modes.}
\label{fig:pipe_err}
\end{figure}

The MAC matrices in Fig.~\ref{fig:pipe_mac} show the correspondence between FOM and reduced modes for the elbow pipe. The estimator performance is illustrated in Fig.~\ref{fig:pipe_err} and Table~\ref{tab:pipe_results}. For the three-substructure tet4 model, the CB error estimator achieves ratios between 0.93 and 1.06, while the HCB-1 estimator achieves ratios between 0.93 and 0.99. The HCB-1 errors are smaller than the CB errors for all 20 modes.

\begin{table}[!htbp]
\centering
\caption{Error estimator validation for the 90-degree elbow pipe (first 10 of 20 modes).}
\label{tab:pipe_results}
\begin{tabular}{cccccccc}
\hline
\multicolumn{2}{c}{} & \multicolumn{3}{c}{CB} & \multicolumn{3}{c}{HCB-1} \\
\cline{3-5}\cline{6-8}
Mode & $f_{\mathrm{FOM}}$ [Hz] & $\xi$ & $\hat{\eta}$ & ratio & $\xi$ & $\hat{\eta}$ & ratio \\
\hline
1  &  200.96 & 1.46e-04 & 1.46e-04 & 0.999 & 3.98e-08 & 3.80e-08 & 0.955 \\
2  &  543.40 & 1.03e-03 & 1.03e-03 & 0.998 & 9.11e-07 & 8.73e-07 & 0.958 \\
3  &  582.41 & 1.44e-03 & 1.44e-03 & 0.997 & 1.21e-06 & 1.17e-06 & 0.966 \\
4  &  668.79 & 2.40e-03 & 2.39e-03 & 0.996 & 1.35e-06 & 1.29e-06 & 0.957 \\
5  &  746.29 & 2.67e-03 & 2.65e-03 & 0.991 & 4.46e-06 & 4.35e-06 & 0.973 \\
6  &  974.38 & 4.68e-03 & 4.75e-03 & 1.015 & 1.72e-06 & 1.64e-06 & 0.953 \\
7  &  999.03 & 4.43e-03 & 4.25e-03 & 0.961 & 3.19e-06 & 3.07e-06 & 0.961 \\
8  & 1014.69 & 4.34e-03 & 4.30e-03 & 0.992 & 5.51e-06 & 5.40e-06 & 0.981 \\
9  & 1138.33 & 6.66e-03 & 6.66e-03 & 0.999 & 1.23e-05 & 1.21e-05 & 0.982 \\
10 & 1173.49 & 2.87e-03 & 2.87e-03 & 1.002 & 2.39e-06 & 2.32e-06 & 0.974 \\
\hline
\end{tabular}
\end{table}

\subsection{Example 3: Reactor pressure vessel}\label{sec:ex_rpv}
\begin{figure}[ht!]
\centering
\includegraphics[width=0.5\textwidth]{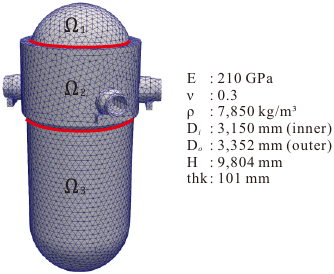}
\caption{Reactor pressure vessel (RPV) example. Geometry including four nozzles (two set-in with bore $\varnothing320$~mm and two set-out) and four support pads, and three-substructure partition along the vessel axis ($y$-direction) at $y = 1{,}500$~mm (Substructure~2/3 boundary) and $y = 3{,}711$~mm (Substructure~1/2 boundary). The tet4 mesh has 26{,}954 elements (8{,}604 nodes, 25{,}782 free DOFs). The upper cylindrical section has a nominal wall thickness of 101~mm. The reported calculations use the coordinate-selected nodal restraints specified in the text.}
\label{fig:rpv_mesh}
\end{figure}

A reactor pressure vessel (RPV) imported from a STEP file is the most complex example. The upper cylindrical section has inner and outer diameters of 3{,}150 and 3{,}352~mm, respectively, corresponding to a wall thickness of 101~mm. The vessel includes four nozzles and four support pads. The STEP geometry is meshed with 26{,}954 tet4 elements and 8{,}604 nodes, as shown in Fig.~\ref{fig:rpv_mesh}. The reported model constrains all three displacement components at 10 nodes selected by a coordinate rule, leaving 25{,}782 free DOFs. The selection takes nodes within 5~mm of the minimum $y$-coordinate among mesh nodes with $2{,}100<y<2{,}500$~mm; the constrained nodes in the stored mesh have $y\approx2{,}100.96$--$2{,}105.65$~mm. The results therefore assess this coordinate-selected nodal restraint configuration and do not represent clamping the complete support-pad bottom faces.

The structure is partitioned into three substructures by element centroid $y$-coordinate (the vessel axis direction), as summarized in Table~\ref{tab:rpv_subs}. The CB reduced-model size is $n_{\mathrm{CB}} = 15 + 15 + 15 + 891 = 936$ (3.6\% of 25{,}782 DOFs).

\begin{table}[!htbp]
\centering
\caption{Substructure configuration for the reactor pressure vessel. Substructure numbering follows Fig.~\ref{fig:rpv_mesh}.}
\label{tab:rpv_subs}
\begin{tabular}{cccc}
\hline
Substructure & Interior DOFs ($N_s^{(k)}$) & Boundary DOFs ($N_b$) & Retained modes ($N_d^{(k)}$) \\
\hline
1 (upper body) & 5{,}541 & -- & 15 \\
2 (middle body) & 8{,}133 & -- & 15 \\
3 (lower body) & 11{,}217 & -- & 15 \\
Boundary & -- & 891 & -- \\
\hline
\end{tabular}
\end{table}

\begin{figure}[ht!]
\centering
\includegraphics[width=\textwidth]{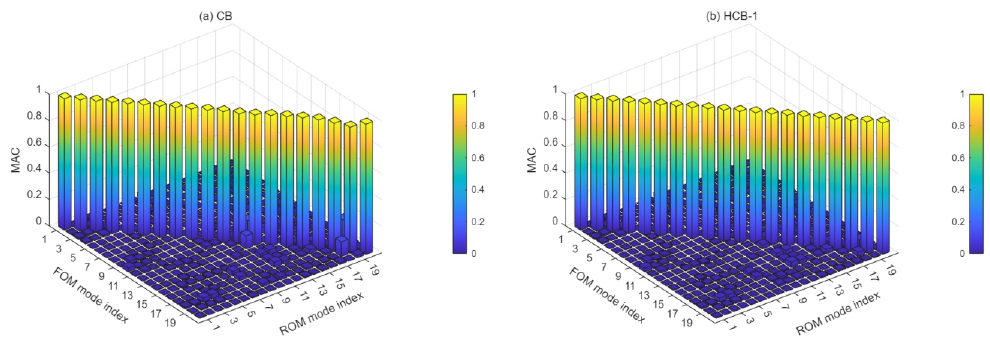}
\caption{Modal assurance criterion $\mathrm{MAC}_{ij}$ for the reactor pressure vessel. (a)~CB and (b)~HCB-1 modes compared with FOM modes for the first 20 modes. MAC is defined in Remark~\ref{sec:error}; values of 1 and 0 indicate aligned and orthogonal modes, respectively.}
\label{fig:rpv_mac}
\end{figure}

\begin{figure}[ht!]
\centering
\includegraphics[width=\textwidth]{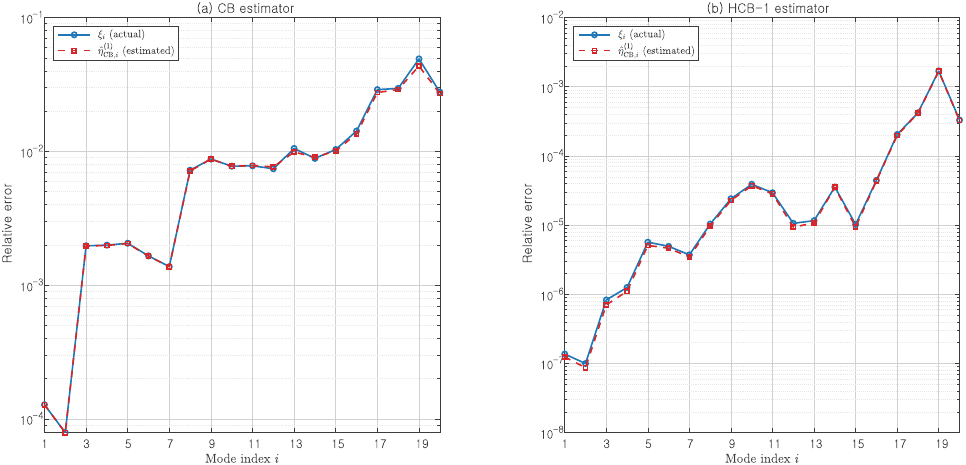}
\caption{Error estimator validation for the reactor pressure vessel. Actual relative eigenvalue error $\xi_i$ (solid blue circles) and estimated error $\hat{\eta}_i$ (dashed red squares) for all 20 modes. (a)~CB estimator, with ratios in the range $0.88$--$1.03$. (b)~HCB-1 estimator, with ratios $0.85$--$0.99$. The HCB-1 errors are smaller than the CB errors for all 20 modes.}
\label{fig:rpv_err}
\end{figure}

The MAC matrices in Fig.~\ref{fig:rpv_mac} show the correspondence between FOM and reduced modes for the RPV. The estimator performance is shown in Fig.~\ref{fig:rpv_err} and Table~\ref{tab:rpv_results}. For the RPV with its complex geometry including nozzles and support pads, the CB error estimator achieves ratios between 0.88 and 1.03, and the HCB-1 estimator achieves ratios between 0.85 and 0.99. The HCB-1 errors are smaller than the CB errors for all 20 modes, showing the accuracy improvement from the residual mode correction for this geometrically complex model.

\begin{table}[ht!]
\centering
\caption{Error estimator validation for the reactor pressure vessel (first 10 of 20 modes).}
\label{tab:rpv_results}
\begin{tabular}{cccccccc}
\hline
\multicolumn{2}{c}{} & \multicolumn{3}{c}{CB} & \multicolumn{3}{c}{HCB-1} \\
\cline{3-5}\cline{6-8}
Mode & $f_{\mathrm{FOM}}$ [Hz] & $\xi$ & $\hat{\eta}$ & ratio & $\xi$ & $\hat{\eta}$ & ratio \\
\hline
1  &  55.62 & 1.28e-04 & 1.28e-04 & 0.999 & 1.37e-07 & 1.24e-07 & 0.908 \\
2  &  56.16 & 7.89e-05 & 7.88e-05 & 0.999 & 9.98e-08 & 8.70e-08 & 0.871 \\
3  & 109.13 & 1.97e-03 & 1.96e-03 & 0.996 & 8.29e-07 & 7.04e-07 & 0.849 \\
4  & 112.24 & 1.99e-03 & 1.98e-03 & 0.995 & 1.26e-06 & 1.11e-06 & 0.881 \\
5  & 132.09 & 2.06e-03 & 2.05e-03 & 0.996 & 5.67e-06 & 5.08e-06 & 0.897 \\
6  & 136.20 & 1.65e-03 & 1.65e-03 & 0.997 & 4.94e-06 & 4.64e-06 & 0.940 \\
7  & 145.04 & 1.38e-03 & 1.37e-03 & 0.999 & 3.72e-06 & 3.45e-06 & 0.926 \\
8  & 154.11 & 7.22e-03 & 7.13e-03 & 0.987 & 1.04e-05 & 9.72e-06 & 0.937 \\
9  & 202.50 & 8.76e-03 & 8.76e-03 & 1.000 & 2.41e-05 & 2.28e-05 & 0.946 \\
10 & 204.56 & 7.75e-03 & 7.73e-03 & 0.997 & 3.89e-05 & 3.71e-05 & 0.953 \\
\hline
\end{tabular}
\end{table}

Table~\ref{tab:cost} reports selected wall-clock times on a desktop workstation (Intel Core i7-13700, 32~GB RAM). The HCB-1 and HCB-2 model timings include construction of the augmented transformation, projection of the stiffness and mass matrices, the augmented eigensolution, SEREP reduction and the subsequent reduced solve and mode expansion. The separate residual-mode row jointly measures the construction of $\mathbf{R}_1$ and $\mathbf{R}_2$. Mesh generation, FE assembly, and partitioning are outside the reported stages.

The recorded subtotal for HCB-2 model construction and solution plus HCB-1 estimator evaluation is approximately $0.13$~s for the plate, $1.99$~s for the pipe, and $17.7$~s for the RPV. These subtotals sum the rounded stage times listed in Table~\ref{tab:cost}.

For construction of both HCB models from the assembled, partitioned system, summing the reported setup rows, the two HCB model rows, and HCB-1 estimator evaluation gives approximately $0.44$~s, $9.96$~s, and $73.7$~s, respectively. This subtotal includes the joint residual-mode construction and excludes the standalone CB eigensolution and CB estimator.

The direct estimator in Eq.~\eqref{eq:eta_direct} requires only the two sets of eigenvalues. The reported HCB-1 estimator evaluation time additionally includes the cross-model mass product, evaluation of $c_i$, and the projection-based diagnostic; it is not a timing of eigenvalue subtraction alone. The CB estimator requires four cross-product matrices ($\mathbf{T}_0^{\mathrm{T}}\mathbf{K}_p\mathbf{T}_r$, $\mathbf{T}_0^{\mathrm{T}}\mathbf{M}_p\mathbf{T}_r$, $\mathbf{T}_r^{\mathrm{T}}\mathbf{K}_p\mathbf{T}_r$, and $\mathbf{T}_r^{\mathrm{T}}\mathbf{M}_p\mathbf{T}_r$), with a reported evaluation cost of 7.4~s for the RPV case.

\begin{table}[!htbp]
\centering
\caption{Reported computational-stage times (seconds) on a desktop workstation (Intel Core i7-13700, 32~GB RAM). HCB model rows include augmented matrix construction and solution plus SEREP, but exclude the joint residual-mode row. Estimator rows assume the required reduced eigensolutions are available. HCB-1 evaluation also includes mode correspondence and the projection diagnostic. Subtotals sum rounded stage times; the plate subtotals omit the evaluation time of less than 0.01~s.}
\label{tab:cost}
\begin{tabular}{lccc}
\hline
Step & Plate & Elbow pipe & RPV \\
\hline
\multicolumn{4}{l}{\textit{CMS setup and basis construction}} \\
\quad Fixed-interface modes & 0.03 & 0.58 & 1.7 \\
\quad Constraint modes & 0.03 & 1.38 & 11.2 \\
\quad Residual modes ($\mathbf{R}_1$, $\mathbf{R}_2$) & 0.06 & 2.78 & 23.2 \\
\quad CB transformation and projected matrices & 0.10 & 1.99 & 10.3 \\
\hline
\multicolumn{4}{l}{\textit{Reduced model construction and solution}} \\
\quad CB eigensolution & 0.02 & 0.04 & 0.10 \\
\quad HCB-1 (augmented model + SEREP) & 0.09 & 1.24 & 9.6 \\
\quad HCB-2 (augmented model + SEREP) & 0.13 & 1.71 & 15.9 \\
\hline
\multicolumn{4}{l}{\textit{Estimator evaluation (eigensolutions available)}} \\
\quad CB estimator (Eq.~\eqref{eq:eta_hcb_efficient}) & 0.02 & 1.18 & 7.4 \\
\quad HCB-1 estimator and mode correspondence & $<$0.01 & 0.28 & 1.8 \\
\hline
\multicolumn{4}{l}{\textit{Reported subtotals (not end-to-end times)}} \\
\quad HCB-2 model + HCB-1 evaluation & 0.13 & 1.99 & 17.7 \\
\quad Setup + both HCB models + evaluation & 0.44 & 9.96 & 73.7 \\
\hline
\end{tabular}
\end{table}

\begin{samepage}
Across the three examples, the HCB-1 estimator provides a lower bound with estimation-to-error ratios of 0.85--0.99, consistent with Section~\ref{sec:morth}. The observed closeness is not guaranteed by the lower-bound property alone. The supporting CB estimator gives ratios of 0.76--1.58 for the plate, 0.93--1.06 for the pipe, and 0.88--1.03 for the RPV over the first 20 modes. Its largest deviation in these examples occurs at plate mode~16, where the ratio is approximately 1.58. The CB estimator uses a leading-order mass approximation and a base-only Rayleigh quotient that need not equal the corresponding CB eigenvalue. These approximations can affect its accuracy; increasing the HCB reference order alone does not guarantee an improvement. Evaluation of the estimators requires no full-order eigensolution.\par
\end{samepage}

\section{Conclusion}\label{sec:conclusion}
This study uses the established HCB enrichment hierarchy to estimate the eigenvalue error remaining in HCB-1. Its contribution is the formulation and analysis of that estimator, beyond the standard CB error targeted by earlier estimators; the HCB reduction method and SEREP construction are inherited from prior work. The main estimate is the relative HCB-1--HCB-2 eigenvalue difference. Under the stated Ritz assumptions, it is a rigorous non-negative lower bound for ordered positive eigenvalues. The approximation $c_i^2\approx1$ is needed only when identifying the related projection-based expression with this direct difference, and plays no role in the lower-bound proof.

The exact error identity separates this bound from a close estimate: the HCB-2 error must be small both in magnitude and relative to the HCB-1 error. The three structural examples give estimation-to-error ratios of approximately 0.85--0.99 for the reported configurations. These observations do not certify a tolerance for an untested model. A decomposition-based extension of the existing CB estimator \cite{kim2014estimating} is a supporting result; its unit-mass and base-only Rayleigh quotient approximations limit its accuracy for some modes.

Neither estimator requires the full-order eigensolution for evaluation. Obtaining an HCB-1 estimate nevertheless requires an HCB-2 reference. The recorded HCB-2 construction/solution and estimator-evaluation subtotals are approximately $0.13$~s, $1.99$~s, and $17.7$~s for the plate, pipe, and RPV. Verification of basis rank and numerical stability and comparison with existing estimators remain necessary for a broader assessment. Closely spaced or repeated eigenvalues also require care in correspondence while preserving the ordered indices used in the bound. The present results concern undamped eigenvalue errors; interface reduction, damping, frequency-response errors, and adaptive tolerance certification are outside the validated scope.

\section*{Data Availability}

Data will be made available on request.

\section*{Acknowledgments}
This work was supported by the National Research Foundation of Korea~(NRF) grant funded by the Korea government~(MSIT) (No.~RS-2026-25501401 and RS-2026-25598696). This research is funded by the ‘Changwon National University - Samsung Changwon Hospital joint Collaboration Research Support Project, South Korea’ in 2025. This study is conducted as part of the Glocal University Project, supported by the RISE (Regional Innovation System \& Education) program funded by the Ministry of Education. This research was supported by the Regional Innovation System \& Education (RISE) program through the RISE Center, Gyeongsangnam-do, funded by the Ministry of Education (MOE) and the Gyeongsangnam-do Provincial Government, Republic of Korea (2026-RISE-16-002). 
% \textbf{Seung-Hwan Boo}: This work was supported by the National Research Foundation of Korea~(NRF) grant funded by the Korea government~(MSIT) (No.~RS-2026-25487009).

\section*{Declaration of Generative AI and AI-assisted technologies in the writing process}

During the preparation of this work, the authors used OpenAI's ChatGPT in order to improve clarity and readability. After using this tool, the authors reviewed and edited the content as needed and take full responsibility for the content of the published article.

\section*{Declaration of competing interest}

The authors declare that they have no known competing financial interests or personal relationships that could have appeared to influence the work reported in this paper.

\section*{CRediT Author Contributions}

\textbf{Jaemin Kim}: Conceptualization, Methodology, Software, Resources, Validation, Formal analysis, Investigation, Data curation, Writing -- original draft, Visualization, Supervision, Project administration, Funding acquisition.
\textbf{Seung-Hwan Boo}: Funding acquisition.

%% ==================================================================
\begin{appendices}
\renewcommand{\theHequation}{\theequation}

\section{Analytical CB reduced matrices}\label{sec:analytical_cb}
\renewcommand{\thefigure}{A\arabic{figure}}
\setcounter{figure}{0}
\renewcommand{\theequation}{A\arabic{equation}}
\setcounter{equation}{0}

The CB reduced matrices can be computed analytically without forming the full transformation product. The reduced stiffness matrix has the block structure
\begin{align}
    \widetilde{\mathbf{K}}_{\mathrm{CB}} = \begin{bmatrix}
    \boldsymbol{\Lambda}_d^{(1)} & & & \mathbf{0} \\
    & \ddots & & \vdots \\
    & & \boldsymbol{\Lambda}_d^{(N_s-1)} & \mathbf{0} \\
    \mathbf{0} & \cdots & \mathbf{0} & \widetilde{\mathbf{K}}_{bb}
    \end{bmatrix},
\end{align}
where the condensed boundary stiffness is
\begin{align}
    \widetilde{\mathbf{K}}_{bb} = \mathbf{K}_{bb} - \sum_{k=1}^{N_s-1} \mathbf{K}_{bs}^{(k)} \left(\mathbf{K}_{ss}^{(k)}\right)^{-1} \mathbf{K}_{sb}^{(k)}.
\end{align}
The reduced mass matrix has the form
\begin{align}
    \widetilde{\mathbf{M}}_{\mathrm{CB}} = \begin{bmatrix}
    \mathbf{I} & & & \widetilde{\mathbf{M}}_{1b} \\
    & \ddots & & \vdots \\
    & & \mathbf{I} & \widetilde{\mathbf{M}}_{(N_s-1)b} \\
    \widetilde{\mathbf{M}}_{b1} & \cdots & \widetilde{\mathbf{M}}_{b(N_s-1)} & \widetilde{\mathbf{M}}_{bb}
    \end{bmatrix},
\end{align}
where the coupling and boundary mass blocks are
\begin{align}
    \widetilde{\mathbf{M}}_{kb} &= (\boldsymbol{\Phi}_d^{(k)})^{\mathrm{T}} \left( \mathbf{M}_{sb}^{(k)} - \mathbf{M}_{ss}^{(k)} \left(\mathbf{K}_{ss}^{(k)}\right)^{-1} \mathbf{K}_{sb}^{(k)} \right), \\
    \widetilde{\mathbf{M}}_{bb} &= \mathbf{M}_{bb} + \sum_{k=1}^{N_s-1} \left[ \mathbf{K}_{bs}^{(k)} \left(\mathbf{K}_{ss}^{(k)}\right)^{-1} \mathbf{M}_{ss}^{(k)} \left(\mathbf{K}_{ss}^{(k)}\right)^{-1} \mathbf{K}_{sb}^{(k)} - \mathbf{M}_{bs}^{(k)} \left(\mathbf{K}_{ss}^{(k)}\right)^{-1} \mathbf{K}_{sb}^{(k)} - \mathbf{K}_{bs}^{(k)} \left(\mathbf{K}_{ss}^{(k)}\right)^{-1} \mathbf{M}_{sb}^{(k)} \right].
\end{align}

\section{Derivation of the residual modal coordinates}\label{sec:derivation_qr}
\renewcommand{\theequation}{B\arabic{equation}}
\setcounter{equation}{0}

This appendix derives Eq.~\eqref{eq:qr_solution} of the manuscript. For brevity, the substructure superscript $(k)$ is dropped throughout.

Substituting the full modal expansion (Eq.~\eqref{eq:full_interior}) into the equation of motion for the interior DOFs (Eq.~\eqref{eq:eom_sub}) gives
\begin{align}\label{eq:D_substituted}
    \mathbf{M}_{ss}\bigl(\boldsymbol{\Phi}_d \ddot{\mathbf{q}}_d + \boldsymbol{\Phi}_r \ddot{\mathbf{q}}_r + \boldsymbol{\Psi}_b \ddot{\mathbf{u}}_b\bigr) + \mathbf{M}_{sb} \ddot{\mathbf{u}}_b + \mathbf{K}_{ss}\bigl(\boldsymbol{\Phi}_d \mathbf{q}_d + \boldsymbol{\Phi}_r \mathbf{q}_r + \boldsymbol{\Psi}_b \mathbf{u}_b\bigr) + \mathbf{K}_{sb} \mathbf{u}_b = \mathbf{0}.
\end{align}
Since $\mathbf{K}_{ss} \boldsymbol{\Psi}_b + \mathbf{K}_{sb} = \mathbf{0}$ by definition of the constraint modes (Eq.~\eqref{eq:constraint}), the stiffness coupling terms cancel, giving
\begin{align}\label{eq:D_simplified}
    \mathbf{M}_{ss}\boldsymbol{\Phi}_d \ddot{\mathbf{q}}_d + \mathbf{M}_{ss}\boldsymbol{\Phi}_r \ddot{\mathbf{q}}_r + \bigl(\mathbf{M}_{ss}\boldsymbol{\Psi}_b + \mathbf{M}_{sb}\bigr) \ddot{\mathbf{u}}_b + \mathbf{K}_{ss}\boldsymbol{\Phi}_d \mathbf{q}_d + \mathbf{K}_{ss}\boldsymbol{\Phi}_r \mathbf{q}_r = \mathbf{0}.
\end{align}
Pre-multiplying by $\boldsymbol{\Phi}_r^{\mathrm{T}}$ and applying the orthogonality conditions---$\boldsymbol{\Phi}_r^{\mathrm{T}} \mathbf{M}_{ss} \boldsymbol{\Phi}_d = \mathbf{0}$, $\boldsymbol{\Phi}_r^{\mathrm{T}} \mathbf{M}_{ss} \boldsymbol{\Phi}_r = \mathbf{I}$, $\boldsymbol{\Phi}_r^{\mathrm{T}} \mathbf{K}_{ss} \boldsymbol{\Phi}_d = \mathbf{0}$, $\boldsymbol{\Phi}_r^{\mathrm{T}} \mathbf{K}_{ss} \boldsymbol{\Phi}_r = \boldsymbol{\Lambda}_r$---yields the decoupled equation for the residual coordinates:
\begin{align}\label{eq:D_decoupled}
    \ddot{\mathbf{q}}_r + \boldsymbol{\Lambda}_r \mathbf{q}_r + \boldsymbol{\Phi}_r^{\mathrm{T}} \hat{\mathbf{M}}_c \ddot{\mathbf{u}}_b = \mathbf{0},
\end{align}
where $\hat{\mathbf{M}}_c = \mathbf{M}_{ss}\boldsymbol{\Psi}_b + \mathbf{M}_{sb}$ is the coupled inertia force (Eq.~\eqref{eq:Mc_hat}). For harmonic motion ($\ddot{\mathbf{q}}_r = -\lambda \mathbf{q}_r$, $\ddot{\mathbf{u}}_b = -\lambda \mathbf{u}_b$), Eq.~\eqref{eq:D_decoupled} becomes
\begin{align}\label{eq:D_harmonic}
    \bigl(\boldsymbol{\Lambda}_r - \lambda \mathbf{I}\bigr) \mathbf{q}_r = \lambda \boldsymbol{\Phi}_r^{\mathrm{T}} \hat{\mathbf{M}}_c \mathbf{u}_b.
\end{align}
Solving for $\mathbf{q}_r$ gives
\begin{align}\label{eq:D_qr_final}
    \mathbf{q}_r = \lambda \bigl(\boldsymbol{\Lambda}_r - \lambda \mathbf{I}\bigr)^{-1} \boldsymbol{\Phi}_r^{\mathrm{T}} \hat{\mathbf{M}}_c \mathbf{u}_b,
\end{align}
which is Eq.~\eqref{eq:qr_solution} of the manuscript.

\section{Derivation of the residual flexibility matrix}\label{sec:derivation_Frs}
\renewcommand{\theequation}{C\arabic{equation}}
\setcounter{equation}{0}

This appendix derives the residual flexibility matrix $\mathbf{F}_1^{(k)}$ used in Eq.~\eqref{eq:res_flex} of the manuscript. For brevity, the substructure superscript $(k)$ is dropped throughout.

The complete set of fixed-interface normal modes satisfies the eigenvalue problem
\begin{align}\label{eq:B_eig_full}
    \mathbf{K}_{ss} \boldsymbol{\phi}_j = \lambda_j \mathbf{M}_{ss} \boldsymbol{\phi}_j, \qquad j = 1, \ldots, N_t,
\end{align}
where $N_t$ is the total number of interior DOFs. Collecting all eigenvectors into $\boldsymbol{\Phi}_{\mathrm{all}} = [\boldsymbol{\phi}_1, \ldots, \boldsymbol{\phi}_{N_t}]$ and the corresponding eigenvalues into $\boldsymbol{\Lambda}_{\mathrm{all}} = \mathrm{diag}(\lambda_1, \ldots, \lambda_{N_t})$, Eq.~\eqref{eq:B_eig_full} is written in matrix form as
\begin{align}\label{eq:B_eig_matrix}
    \mathbf{K}_{ss} \boldsymbol{\Phi}_{\mathrm{all}} = \mathbf{M}_{ss} \boldsymbol{\Phi}_{\mathrm{all}} \boldsymbol{\Lambda}_{\mathrm{all}}.
\end{align}
Since the $N_t$ eigenvectors form a complete set, $\boldsymbol{\Phi}_{\mathrm{all}}$ is square and invertible. With mass normalization
\begin{align}\label{eq:B_M_ortho}
    \boldsymbol{\Phi}_{\mathrm{all}}^{\mathrm{T}} \mathbf{M}_{ss} \boldsymbol{\Phi}_{\mathrm{all}} = \mathbf{I},
\end{align}
the inverse is identified as $\boldsymbol{\Phi}_{\mathrm{all}}^{-1} = \boldsymbol{\Phi}_{\mathrm{all}}^{\mathrm{T}} \mathbf{M}_{ss}$, which also implies the resolution of identity $\boldsymbol{\Phi}_{\mathrm{all}} \boldsymbol{\Phi}_{\mathrm{all}}^{\mathrm{T}} = \mathbf{M}_{ss}^{-1}$. Right-multiplying Eq.~\eqref{eq:B_eig_matrix} by $\boldsymbol{\Phi}_{\mathrm{all}}^{-1} = \boldsymbol{\Phi}_{\mathrm{all}}^{\mathrm{T}} \mathbf{M}_{ss}$ yields the spectral decomposition of the stiffness matrix:
\begin{align}\label{eq:B_spectral_K}
    \mathbf{K}_{ss} = \mathbf{M}_{ss} \boldsymbol{\Phi}_{\mathrm{all}} \boldsymbol{\Lambda}_{\mathrm{all}} \boldsymbol{\Phi}_{\mathrm{all}}^{\mathrm{T}} \mathbf{M}_{ss}.
\end{align}
Similarly, pre-multiplying Eq.~\eqref{eq:B_eig_matrix} by $\boldsymbol{\Phi}_{\mathrm{all}}^{\mathrm{T}}$ and applying the mass orthonormality~\eqref{eq:B_M_ortho} yields the stiffness orthogonality:
\begin{align}\label{eq:B_K_ortho}
    \boldsymbol{\Phi}_{\mathrm{all}}^{\mathrm{T}} \mathbf{K}_{ss} \boldsymbol{\Phi}_{\mathrm{all}} = \boldsymbol{\Lambda}_{\mathrm{all}}.
\end{align}
Inverting both sides of Eq.~\eqref{eq:B_K_ortho}, i.e.\ $\boldsymbol{\Phi}_{\mathrm{all}}^{-1} \mathbf{K}_{ss}^{-1} (\boldsymbol{\Phi}_{\mathrm{all}}^{\mathrm{T}})^{-1} = \boldsymbol{\Lambda}_{\mathrm{all}}^{-1}$, and pre-multiplying by $\boldsymbol{\Phi}_{\mathrm{all}}$ and post-multiplying by $\boldsymbol{\Phi}_{\mathrm{all}}^{\mathrm{T}}$ gives the full flexibility matrix:
\begin{align}\label{eq:B_full_flex}
    \mathbf{K}_{ss}^{-1} = \boldsymbol{\Phi}_{\mathrm{all}} \boldsymbol{\Lambda}_{\mathrm{all}}^{-1} \boldsymbol{\Phi}_{\mathrm{all}}^{\mathrm{T}} = \sum_{j=1}^{N_t} \frac{1}{\lambda_j} \boldsymbol{\phi}_j \boldsymbol{\phi}_j^{\mathrm{T}}.
\end{align}
In general, $\boldsymbol{\Phi}_{\mathrm{all}}^{-1} = \boldsymbol{\Phi}_{\mathrm{all}}^{\mathrm{T}} \mathbf{M}_{ss} \neq \boldsymbol{\Phi}_{\mathrm{all}}^{\mathrm{T}}$, since $\mathbf{M}_{ss} \neq \mathbf{I}$. The appearance of the plain transpose $\boldsymbol{\Phi}_{\mathrm{all}}^{\mathrm{T}}$ (without $\mathbf{M}_{ss}$) in Eq.~\eqref{eq:B_full_flex} is a direct consequence of the stiffness orthogonality~\eqref{eq:B_K_ortho}, not of standard orthogonality.

The modes are partitioned into $N_d$ dominant modes and $N_r = N_t - N_d$ residual (truncated) modes:
\begin{align}\label{eq:B_partition_modes}
    \boldsymbol{\Phi}_{\mathrm{all}} = [\underbrace{\boldsymbol{\Phi}_d}_{N_d\;\mathrm{dominant}},\; \underbrace{\boldsymbol{\Phi}_r}_{N_r\;\mathrm{residual}}], \qquad
    \boldsymbol{\Lambda}_{\mathrm{all}} = \begin{bmatrix} \boldsymbol{\Lambda}_d & \\ & \boldsymbol{\Lambda}_r \end{bmatrix}.
\end{align}
Substituting this partition into Eq.~\eqref{eq:B_full_flex} yields
\begin{align}\label{eq:B_flex_split}
    \mathbf{K}_{ss}^{-1}
    = \underbrace{\boldsymbol{\Phi}_d \boldsymbol{\Lambda}_d^{-1} \boldsymbol{\Phi}_d^{\mathrm{T}}}_{\text{dominant flexibility}}
    + \underbrace{\boldsymbol{\Phi}_r \boldsymbol{\Lambda}_r^{-1} \boldsymbol{\Phi}_r^{\mathrm{T}}}_{\text{residual flexibility}}.
\end{align}

Rearranging Eq.~\eqref{eq:B_flex_split}, the residual flexibility matrix is obtained as
\begin{align}\label{eq:B_Frs_def}
    \mathbf{F}_1 = \mathbf{K}_{ss}^{-1} - \boldsymbol{\Phi}_d \boldsymbol{\Lambda}_d^{-1} \boldsymbol{\Phi}_d^{\mathrm{T}} = \boldsymbol{\Phi}_r \boldsymbol{\Lambda}_r^{-1} \boldsymbol{\Phi}_r^{\mathrm{T}} = \sum_{j=N_d+1}^{N_t} \frac{1}{\lambda_j} \boldsymbol{\phi}_j \boldsymbol{\phi}_j^{\mathrm{T}}.
\end{align}
This matrix represents the flexibility contribution of the truncated modes only. It can be computed from $\mathbf{K}_{ss}^{-1}$ and the dominant eigenpairs $(\boldsymbol{\Phi}_d, \boldsymbol{\Lambda}_d)$ without explicitly computing the residual eigenpairs $(\boldsymbol{\Phi}_r, \boldsymbol{\Lambda}_r)$.

To capture higher-order inertia effects, the second-order residual flexibility is obtained by applying the mass operator between two copies of $\mathbf{F}_1$. Using the spectral form in Eq.~\eqref{eq:B_Frs_def} and the mass orthogonality $\boldsymbol{\Phi}_r^{\mathrm{T}} \mathbf{M}_{ss} \boldsymbol{\Phi}_r = \mathbf{I}$ gives
\begin{align}\label{eq:B_Frm_derivation}
    \mathbf{F}_2
    &= \mathbf{F}_1 \mathbf{M}_{ss} \mathbf{F}_1 \notag \\
    &= \left(\boldsymbol{\Phi}_r \boldsymbol{\Lambda}_r^{-1} \boldsymbol{\Phi}_r^{\mathrm{T}}\right) \mathbf{M}_{ss} \left(\boldsymbol{\Phi}_r \boldsymbol{\Lambda}_r^{-1} \boldsymbol{\Phi}_r^{\mathrm{T}}\right) \notag \\
    &= \boldsymbol{\Phi}_r \boldsymbol{\Lambda}_r^{-1} \underbrace{\left(\boldsymbol{\Phi}_r^{\mathrm{T}} \mathbf{M}_{ss} \boldsymbol{\Phi}_r\right)}_{=\mathbf{I}} \boldsymbol{\Lambda}_r^{-1} \boldsymbol{\Phi}_r^{\mathrm{T}} \notag \\
    &= \boldsymbol{\Phi}_r \boldsymbol{\Lambda}_r^{-2} \boldsymbol{\Phi}_r^{\mathrm{T}} = \sum_{j=N_d+1}^{N_t} \frac{1}{\lambda_j^2} \boldsymbol{\phi}_j \boldsymbol{\phi}_j^{\mathrm{T}}.
\end{align}
Similarly, this can be computed without the residual eigenpairs using:
\begin{align}\label{eq:B_Frm_practical}
    \mathbf{F}_2 = \mathbf{F}_1 \mathbf{M}_{ss} \mathbf{F}_1 = \mathbf{K}_{ss}^{-1} \mathbf{M}_{ss} \mathbf{K}_{ss}^{-1} - \boldsymbol{\Phi}_d \boldsymbol{\Lambda}_d^{-2} \boldsymbol{\Phi}_d^{\mathrm{T}}.
\end{align}
This pattern generalizes: the $n$-th order residual flexibility is $\mathbf{F}_{r}^{(n)} = \boldsymbol{\Phi}_r \boldsymbol{\Lambda}_r^{-n} \boldsymbol{\Phi}_r^{\mathrm{T}}$, with each successive order weighting lower-frequency residual modes more heavily, thereby improving the representation of the inertia effects neglected by mode truncation.

\section{Derivation of the decomposition-based error estimator}\label{sec:derivation_decomp}
\renewcommand{\theequation}{D\arabic{equation}}
\setcounter{equation}{0}

This appendix derives Eq.~\eqref{eq:eta_enhanced} of the manuscript from the decomposition $\mathbf{u}_i = \mathbf{u}_{0,i} + \mathbf{u}_{r,i}$ (Eq.~\eqref{eq:decomp_general}).
Let $(\bar{\lambda}_i, \bar{\boldsymbol{\phi}}_i)$ denote the $i$-th eigenpair of the enhanced (HCB) reduced system (Eq.~\eqref{eq:serep_reduced}), so that the corresponding full-order approximation is $\mathbf{u}_i = \mathbf{T}_{\mathrm{comp}} \bar{\boldsymbol{\phi}}_i$ with mass-normalization $\mathbf{u}_i^{\mathrm{T}} \mathbf{M}_p \mathbf{u}_i = 1$.
The base-only component $\mathbf{u}_{0,i} = \mathbf{T}_0 \bar{\boldsymbol{\phi}}_i$ and the residual component $\mathbf{u}_{r,i} = \mathbf{T}_r \bar{\boldsymbol{\phi}}_i$ (Eq.~\eqref{eq:decomp_general}) satisfy $\mathbf{u}_i = \mathbf{u}_{0,i} + \mathbf{u}_{r,i}$.
Define $\hat{\lambda}_{0,i} \coloneqq \mathbf{u}_{0,i}^{\mathrm{T}} \mathbf{K}_p \mathbf{u}_{0,i} / \mathbf{u}_{0,i}^{\mathrm{T}} \mathbf{M}_p \mathbf{u}_{0,i}$ as the Rayleigh quotient of the base-only approximation.

From the mass-normalization and the eigenvalue relation, $\mathbf{u}_i^{\mathrm{T}} \mathbf{K}_p \mathbf{u}_i = \bar{\lambda}_i$. Writing $\mathbf{u}_{0,i} = \mathbf{u}_i - \mathbf{u}_{r,i}$ and expanding the stiffness and mass quadratic forms gives
\begin{align}
    \mathbf{u}_{0,i}^{\mathrm{T}} \mathbf{K}_p \mathbf{u}_{0,i} &= \bar{\lambda}_i - 2\mathbf{u}_{r,i}^{\mathrm{T}} \mathbf{K}_p \mathbf{u}_i + \mathbf{u}_{r,i}^{\mathrm{T}} \mathbf{K}_p \mathbf{u}_{r,i}, \label{eq:E_K_expand} \\
    \mathbf{u}_{0,i}^{\mathrm{T}} \mathbf{M}_p \mathbf{u}_{0,i} &= 1 - 2\mathbf{u}_{r,i}^{\mathrm{T}} \mathbf{M}_p \mathbf{u}_i + \mathbf{u}_{r,i}^{\mathrm{T}} \mathbf{M}_p \mathbf{u}_{r,i}. \label{eq:E_M_expand}
\end{align}
The relative difference between the base-only Rayleigh quotient and the enhanced eigenvalue is
\begin{align}\label{eq:E_rel_error}
    \delta_{0,i} \coloneqq \frac{\hat{\lambda}_{0,i} - \bar{\lambda}_i}{\bar{\lambda}_i}
    = \frac{\mathbf{u}_{0,i}^{\mathrm{T}} \mathbf{K}_p \mathbf{u}_{0,i} - \bar{\lambda}_i \, \mathbf{u}_{0,i}^{\mathrm{T}} \mathbf{M}_p \mathbf{u}_{0,i}}{\bar{\lambda}_i \, \mathbf{u}_{0,i}^{\mathrm{T}} \mathbf{M}_p \mathbf{u}_{0,i}}.
\end{align}
Substituting Eqs.~\eqref{eq:E_K_expand}--\eqref{eq:E_M_expand} and simplifying using $\mathbf{u}_i = \mathbf{u}_{0,i} + \mathbf{u}_{r,i}$, the numerator reduces to
\begin{align}
    -\bigl[2\mathbf{u}_{0,i}^{\mathrm{T}} (\mathbf{K}_p - \bar{\lambda}_i \mathbf{M}_p) \mathbf{u}_{r,i} + \mathbf{u}_{r,i}^{\mathrm{T}} (\mathbf{K}_p - \bar{\lambda}_i \mathbf{M}_p) \mathbf{u}_{r,i}\bigr]. \label{eq:E_numerator}
\end{align}
Dividing by $\bar{\lambda}_i \, \mathbf{u}_{0,i}^{\mathrm{T}} \mathbf{M}_p \mathbf{u}_{0,i}$ gives the exact expression
\begin{align}\label{eq:E_eta_raw}
    \delta_{0,i}
    = -\frac{2\mathbf{u}_{0,i}^{\mathrm{T}} (\mathbf{K}_p - \bar{\lambda}_i \mathbf{M}_p) \mathbf{u}_{r,i} + \mathbf{u}_{r,i}^{\mathrm{T}} (\mathbf{K}_p - \bar{\lambda}_i \mathbf{M}_p) \mathbf{u}_{r,i}}{\bar{\lambda}_i \, \mathbf{u}_{0,i}^{\mathrm{T}} \mathbf{M}_p \mathbf{u}_{0,i}}.
\end{align}
Since $\mathbf{u}_i^{\mathrm{T}} \mathbf{M}_p \mathbf{u}_i = 1$ and $\mathbf{u}_{r,i}$ is small, $\mathbf{u}_{0,i}^{\mathrm{T}} \mathbf{M}_p \mathbf{u}_{0,i} \approx 1$ to leading order. Using the auxiliary matrix $\mathbf{G}(\mu) = \mathbf{M}_p - \mu^{-1}\mathbf{K}_p$ (Eq.~\eqref{eq:Gi}), note that
\begin{align}
    -\frac{1}{\bar{\lambda}_i}(\mathbf{K}_p - \bar{\lambda}_i \mathbf{M}_p) = \mathbf{M}_p - \frac{1}{\bar{\lambda}_i}\mathbf{K}_p = \mathbf{G}(\bar{\lambda}_i). \notag
\end{align}
Applying this identity to Eq.~\eqref{eq:E_eta_raw} with the leading-order approximation $\mathbf{u}_{0,i}^{\mathrm{T}} \mathbf{M}_p \mathbf{u}_{0,i} \approx 1$ gives
\begin{align}\label{eq:E_eta_simplified}
    \delta_{0,i}
    \approx 2\mathbf{u}_{0,i}^{\mathrm{T}} \mathbf{G}(\bar{\lambda}_i) \, \mathbf{u}_{r,i} + \mathbf{u}_{r,i}^{\mathrm{T}} \mathbf{G}(\bar{\lambda}_i) \, \mathbf{u}_{r,i}.
\end{align}
Defining the leading-order estimator by the right-hand side and substituting $\mathbf{u}_{0,i} = \mathbf{T}_0 \bar{\boldsymbol{\phi}}_i$ and $\mathbf{u}_{r,i} = \mathbf{T}_r \bar{\boldsymbol{\phi}}_i$ gives Eq.~\eqref{eq:eta_enhanced}:
\begin{align}\label{eq:E_eta_final}
    \hat{\eta}_i = 2\bar{\boldsymbol{\phi}}_i^{\mathrm{T}} \mathbf{T}_0^{\mathrm{T}} \mathbf{G}(\bar{\lambda}_i) \, \mathbf{T}_r \bar{\boldsymbol{\phi}}_i + \bar{\boldsymbol{\phi}}_i^{\mathrm{T}} \mathbf{T}_r^{\mathrm{T}} \mathbf{G}(\bar{\lambda}_i) \, \mathbf{T}_r \bar{\boldsymbol{\phi}}_i.
\end{align}

\end{appendices}

\bibliographystyle{elsarticle-num}
\bibliography{references}

\end{document}